\documentclass[10pt]{ijnam}
\def \qed {\hfill \rule{1ex}{1ex}}
\def \U {{\bf u}}
\def \x {{\bf x}}

\def \H {{\bf H}}

\def \V {{\bf v}}
\def \bt {{\widetilde{\hbox{\rm \bf b}}}}
\def \W {{\bf w}}
\def \N {{\mathbf{n}}}

\def \Ut {{\tilde {\mathbf{u}}}}
\def \Deltat {{\widetilde \Delta}}
\def \F {{\bf f }}
\def \P {{\bf P}}
\def \R {{\mathbf{R}}}

\def \L {{\bf L}}

\def \b {{\hbox{\rm  b}}}

\def \Deltat{{{\mathbf{\widetilde \Delta}}}}
\def \Nabla {{\hbox{\boldmath $\nabla$ \unboldmath \!\!}}}
\def \Del {{\hbox{\boldmath $\delta$ \unboldmath \!\!}}}
\def \psig {{\hbox{\boldmath $\psi$ \unboldmath \!\!}}}

\newtheorem{lem}{Lemma}[section]
\newtheorem{theo}{Theorem}[section]
\newtheorem{prop}{Proposition}[section]

\makeatletter
\renewcommand\theequation{\thesection.\arabic{equation}}
\@addtoreset{equation}{section} \makeatother

\def\currentvolume{1}
\def\currentissue{1}
\def\currentyear{2007}
\def\month{April}
\pagespan{1}{18}
\copyrightinfo{2007}{} 

\begin{document}

\title[A colocated finite volume scheme for the Navier-Stokes equations]
{Stability of a colocated finite volume scheme\\
for the incompressible Navier-Stokes equations}

\author[S. Zimmermann]{S\'ebastien Zimmermann}
\address{
  Department of Mathematics,
  Centrale Lyon University,
  63177 Ecully, FRANCE
} \email{Sebastien.Zimmermann@ec-lyon.fr}



\commby{Jean-Luc Guermond}

\date{April 1, 2007 and, in revised form, April 1, 2007.}



\abstract{We introduce a finite volume scheme for the
two-dimensional incompressible Navier-Stokes equations. We use a
triangular mesh. The unknowns for the velocity and pressure are both
piecewise constant (colocated scheme). We use a projection
(fractional-step) method to deal with the incompressibility
constraint. We prove that the differential operators in the
Navier-Stokes equations and their discrete counterparts share similar
properties. In particular, we state an inf-sup (Babu\v{s}ka-Brezzi)
condition. We infer from it the stability of the scheme.}


\maketitle

\section{Introduction}

We consider the flow of an incompressible fluid in a open bounded
set $\Omega \subset \mathbb{R}^2$ during the time interval $[0,T]$.
The velocity field $\mathbf{u}:\Omega \times [0,T] \to \mathbb{R}^2
$ and the pressure field $p:\Omega \times [0,T] \to \mathbb{R}$
satisfy the Navier-Stokes equations
\begin{eqnarray}
&& \mathbf{u}_t - \frac{1}{\hbox{Re}} \, \mathbf{\Delta} \mathbf{u}
+(\U \cdot \Nabla) \U+\nabla p = \F \, ,
\label{eq:mom} \\
&& \hbox{div } \U=0 \, , \label{eq:incomp}
\end{eqnarray}
with the boundary and initial condition
\begin{equation*}
 \U|_{\partial \Omega}=0 \, , \hspace{2cm}  \U|_{t=0}=\U_0.
\end{equation*}
The terms $\mathbf{\Delta} \U$ and $(\U \cdot \Nabla)\U$ are
respectively associated with the  physical phenomena of diffusion
and convection. The Reynolds number $\hbox{Re}$ measures the
influence of convection in the flow. For equations
(\ref{eq:mom})--(\ref{eq:incomp}), finite element and finite
difference methods are well known and mathematical studies are
available (see \cite{giraultr} for example). Numerous computations
have also been conducted with finite volume schemes (\emph{e.g.}
\cite{kimchoi} and \cite{boicaya}).
 However, in this case, few mathematical results are available.  Let us cite {\sc Eymard and Herbin} \cite{herb3}
 and {\sc Eymard, Latch\'e and Herbin} \cite{herb4}.
 In order to deal with the incompressibility constraint (\ref{eq:incomp}), these works use a penalization  method.
 Another way is to use  the projection methods which have been
 introduced by {\sc Chorin} \cite{chorin} and {\sc Temam} \cite{temam}.
 This is the case in  {\sc Faure} \cite{faure}. In this work, however,
 the mesh is made of  squares, so that  the geometry of the problem is limited.
 Therefore, we introduce in what follows a finite volume scheme on triangular meshes for equations
(\ref{eq:mom})--(\ref{eq:incomp}), using a projection method. An
interesting feature of this scheme is that the unknowns for the
velocity and pressure are both piecewise constant (colocated
scheme). It leads to an economic computer storage, and allows an
easy generalization of the scheme to the 3D case. The layout of the
article is the following. We first introduce (section
\ref{sec:notd}) some notations and hypotheses on the mesh. We define
(section \ref{subsec:espd}) the spaces we use to approximate the
velocity and pressure. We define also (section \ref{subsec:opd}) the
operators we use to approximate the differential operators in
(\ref{eq:mom})--(\ref{eq:incomp}).
 Combining this with a projection method, we build the scheme in section \ref{sec:presschema}.
  In order to provide a mathematical analysis for the scheme, we
  prove in section \ref{sec:propopd} that the differential operators
   in (\ref{eq:mom})--(\ref{eq:incomp}) and their discrete counterparts share similar properties.
 In particular, the discrete operators for the gradient and the
 divergence are adjoint. Also, the discrete gradient operator is a consistent approximation of its continuous counterpart.
  The discrete operator for the convection term
 is positive, stable and consistent. The discrete operator for the divergence
  satisfies an inf-sup (Babu\v{s}ka-Brezzi) condition.
From these properties we deduce in section \ref{sec:stab} the
stability of the scheme.

\noindent We conclude with some notations. The spaces $(L^2,|.|)$ and
$(L^\infty,\|.\|_\infty)$ are the usual Lebesgue spaces and we set
$L^2_0=\{q \in L^2 \, ; \int_\Omega q(\x) \, d\x=0 \}$. Their
vectorial counterparts are $(\L^2,|.|)$ and
$(\L^\infty,\|.\|_\infty)$ with $\L^2=(L^2)^2$ and
$\L^\infty=(L^\infty)^2$. For $k \in \mathbb{N}^*$,
$(H^k,\|\cdot\|_k)$ is the usual Sobolev space. Its vectorial
counterpart is $(\H^k,\|.\|_k)$ with $\H^k=(H^k)^2$. For $k=1$, the
functions of $\H^1$ with a null trace on the boundary form the space
$\H^1_0$. Also, we set $\Nabla \U=(\nabla u_1,\nabla u_2)^T$ if
$\U=(u_1,u_2) \in \H^1$.
If $\mathbf{X}\subset \L^2$ is a Banach space, we define
${\mathcal{C}}(0,T;\mathbf{X})$
 (resp. $L^2(0,T;\mathbf{X})$) as the set of the
 applications $\mathbf{g}:[0,T] \to \mathbf{X}$ such that
$t \to |\mathbf{g}(t)|$ is continous (resp. square integrable).
 The norms $\|.\|_{{\mathcal{C}}(0,T;\mathbf{X})}$ and
$\|.\|_{L^2(0,T;\mathbf{X})}$
are defined respectively by
$\|\mathbf{g}\|_{{\mathcal{C}}(0,T;\mathbf{X})}=\sup_{t\in[0,T]}
|\mathbf{g}(t)|$ and $\|\mathbf{g}\|_{L^2(0,T;\mathbf{X})}=\left(
\int^T_0 |g(t)|^2 \, ds \right)^{1/2}$ .
In all calculations, $C$ is a generic positive constant, depending
only on $\Omega$, $\U_0$ and $\F$.

\section{Discrete setting}
\label{sec:notd}

First, we introduce the spaces and the operators needed to build the
scheme.

\subsection{The mesh}

Let ${\mathcal{T}}_h$ be a triangular mesh of $\Omega$:
$\overline{\Omega}=\cup_{K \in {\mathcal{T}}_h} K $. For each
triangle $K \in {\mathcal{T}}_h$, we denote by  $|K|$ its area and
${\mathcal{E}}_K$ the set of his edges.
 If $\sigma \in {\mathcal{E}}_K$, $\N_{K,\sigma}$  is the unit vector normal to $\sigma$ pointing outward of $K$.

\noindent  The set of edges of the mesh is ${\mathcal{E}}_h=\cup_{K
\in {\mathcal{T}}_h} {\mathcal{E}}_K$. The length of an edge $\sigma
\in {\mathcal{E}}_h$ is $|\sigma|$ and its middle point $\x_\sigma$.
The set of edges located inside $\Omega$ (resp. on its boundary) is
${\mathcal{E}}^{int}_h$
 (resp. ${\mathcal{E}}^{ext}_h$): ${\mathcal{E}}_h={\mathcal{E}}^{int}_h \cup {\mathcal{E}}^{ext}_h$.
 If $\sigma \in {\mathcal{E}}^{int}_h$, $K_\sigma$ and $L_\sigma$
 are the triangles sharing $\sigma$ as an edge.
If $\sigma \in {\mathcal{E}}^{ext}_h$, only the triangle $K_\sigma$
 inside $\Omega$ is defined.

\noindent We denote by   $\x_K$  the circumcenter of a triangle $K$.
We assume that the measure of all interior angles of the triangles
of the mesh are below $\frac{\pi}{2}$, so that $\x_K \in K$. If
$\sigma \in {\mathcal{E}}^{int}_h$ (resp. $\sigma \in
{\mathcal{E}}^{ext}_h$ ) we set
 $d_\sigma=d(\x_{K_\sigma},\x_{L_\sigma})$ (resp.  $d_\sigma=d(\x_\sigma,\x_{K_\sigma})$).
We define for all edge $\sigma \in {\mathcal{E}}_h$
\begin{equation}
\label{eq:deftaus}
  \tau_\sigma=\frac{|\sigma|}{d_\sigma}.
\end{equation}
The maximum circumradius  of the triangles of the mesh is $h$. We assume  (\cite{eymgal} p. 776) that there exists $C>0$  such that
\begin{equation*}
 \forall \, \sigma \in {\mathcal{E}}_h, \hspace{1cm} \; \; d(\x_{K_\sigma},\sigma) \ge C |\sigma|
\hspace{.5cm}   \hbox{and} \hspace{.5cm} |\sigma| \ge Ch.
\end{equation*}
It implies that  there exists $C>0$ such that
\begin{equation}
\label{eq:proptaus}
  \forall \, \sigma \in {\mathcal{E}}_h \,, \hspace{1cm} \tau_\sigma \ge
  C \, ,
\end{equation}
and  for all triangles  $K \in {\mathcal{T}}_h$ we have (with
$\sigma \in {\mathcal{E}}_K$ and $h_{K,\sigma}$ the matching
altitude)
\begin{equation}
\label{eq:propairet}
  |K|=\frac{1}{2} \, |\sigma| \, h_{K,\sigma} \ge \frac{1}{2} \, |\sigma| \, d(\x_K,\x_\sigma) \ge C \, h^2.
\end{equation}
Lastly, if $K \in {\mathcal{T}}_h$ and $L \in {\mathcal{T}}_h$ are
two triangles sharing the edge
 $\sigma \in {\mathcal{E}}^{int}_h$, we define
 \begin{equation*}
  \alpha_{K,L}=\frac{d(\x_L,\x_\sigma)}{d(\x_K,\x_L)}.
 \end{equation*}
Let us notice that $\alpha_{K,L} \in [0,1]$ and
$\alpha_{K,L}+\alpha_{L,K}=1$.

\subsection{The discrete spaces}

\label{subsec:espd}

\noindent   We first define
 \begin{equation*}
 \label{eq:defp0}
   P_0= \{ q \in L^2 \; ; \; \forall \, K \in {\mathcal{T}}_h, \; \; q|_K \hbox{ is a constant} \} \, ,
   \hspace{1cm} \P_0=(P_0)^2.
\end{equation*}
For the sake of concision, we set for all $q_h \in P_0$ (resp. $\V_h
\in \P_0$) and all triangle $K \in {\mathcal{T}}_h$:
 $q_K=q_h|_K$ (resp. $\V_K=\V_h|_K$).
Although $\P_0 \not \subset \H^1$, we define the discrete equivalent
of a $\H^1$ norm as follows. For all $\V_h \in \P_0$ we set
\begin{equation}
\label{eq:defh1d}
  \|\V_h\|_h=\left( \sum_{\sigma \in {\mathcal{E}}^{int}_h} \tau_{\sigma} \, |\V_{L_\sigma} - \V_{K_\sigma}|^2
  +\sum_{\sigma \in {\mathcal{E}}^{ext}_h} \tau_\sigma \, |\V_{K_\sigma}|^2 \right)^{1/2}
\end{equation}
where $\tau_\sigma$ is given by (\ref{eq:deftaus}).
 We have \cite{eymgal} a Poincar\'e-like inequality for $\P_0$: there exists $C>0$ such that for all
  $\V_h \in \P_0$
 \begin{equation}
 \label{eq:inpoinp0}
  |\V_h| \le C \, \|\V_h\|_h.
 \end{equation}
We also have the following inverse inequality.
\begin{prop}
\label{propinvp0}

\noindent There exists a constant $C>0$ such that for all $\V_h \in
\P_0$
\begin{equation*}
\label{eq:iinv} h \, \|\V_h\|_h \le C \, |\V_h|.
\end{equation*}
\end{prop}
\noindent {\sc Proof.} According to (\ref{eq:defh1d})
\begin{equation*}
h^2 \,  \|\V_h\|^2_h = \sum_{\sigma \in {\mathcal{E}}^{int}_h} h^2
\, \tau_\sigma \, |\V_{L_\sigma} - \V_{K_\sigma}|^2 +\sum_{\sigma
\in {\mathcal{E}}^{ext}_h} h^2 \, \tau_\sigma \, |\V_{K_\sigma}|^2.
\end{equation*}
We deduce from (\ref{eq:proptaus}) and (\ref{eq:propairet}) that
$  h^2  \, \tau_\sigma
  \le C  \, |K_\sigma|$ and  $h^2  \, \tau_\sigma  \le C  \,
  |L_\sigma|$.
Thus, since
$|\V_{L_\sigma} - \V_{K_\sigma}|^2 \le 2 \, \big( |\V_{L_\sigma}|^2
+|\V_{K_\sigma}|^2\big)$,
we get
\begin{equation*}
 h^2 \, \|\V_h\|^2_h \le C\sum_{\sigma \in {\mathcal{E}}^{int}_h} \big(|K_\sigma| \, |\V_{K_\sigma}|^2
+ |L_\sigma| \, |\V_{L_\sigma}|^2\big) +C\sum_{\sigma \in
{\mathcal{E}}^{ext}_h} |K_\sigma| \, |\V_{K_\sigma}|^2.
\end{equation*}
Hence
$ h^2 \, \|\V_h\|^2_h   \le  C\sum_{K \in {\mathcal{T}}_h} |K| \,
|\V_K|^2
 \le  C \, |\V_h|^2$.  \qed
\vspace{.1cm}

\noindent From the norm $\|.\|_h$ we deduce a dual norm.
 For all $\V_h \in \P_0$ we set
\begin{equation}
\label{eq:defnormdp0}
  \|\V_h\|_{-1,h}=\sup_{\psig_h \in \P_0} \frac{(\V_h,\psig_h)}{\|\psig_h\|_h}.
\end{equation}
For all $\U_h \in \P_0$ and $\V_h \in \P_0$ we have
  $(\U_h,\V_h) \le \|\U_h\|_{-1,h} \, \|\V_h\|_h$.
Now we  introduce some operators on $P_0$ and $\P_0$. We define the
projection operator $\Pi_{\P_0}: \L^2 \to \P_0$ as follows.
 For all $\W \in \L^2$, $\Pi_{\P_0} \W \in \P_0$ is given by
\begin{equation}
\label{eq:defprojp0}
   \forall \, K \in {\mathcal{T}}_h \, ,  \hspace{1cm} (\Pi_{\P_0} \W)|_K=\frac{1}{|K|} \int_K \W(\x) \, d\x.
\end{equation}
We easily check that for all  $\W \in \L^2$ and $\V_h \in \P_0$ we
have $(\Pi_{\P_0} \W, \V_h)=(\W,\V_h)$. It implies that $\Pi_{\P_0}$
is stable for the $\L^2$ norm.
 We define also the interpolation operator $\widetilde \Pi_{P_0}:H^2 \to P_0$.
 For all $q\in H^2$, $\widetilde \Pi_{P_0} q \in P_0$
 is given by
\begin{equation*}
  \forall \, K \in {\mathcal{T}}_h \, , \hspace{1cm} \widetilde \Pi_{P_0} q|_K=q(\x_K).
\end{equation*}
According to the Sobolev embedding theorem, $q \in H^2$ is a.e.
equal to a continuous function. Therefore the definition above makes
sense. We also set $\widetilde \Pi_{\P_0}=(\widetilde \Pi_{P_0})^2$.
The operator $\widetilde \Pi_{P_0}$ (resp. $\widetilde \Pi_{\P_0}$)
is naturally stable for the $L^\infty$ (resp. $\L^\infty$) norm. One also
checks  (\cite{brenns} and \cite{zimm}) that there exists $C>0$ such that
\begin{equation}
\label{eq:estdiffprojp0}
 |\V-\Pi_{\P_0} \V| \le C\, h \, \|\V\|_1 \,, \hspace{1cm}|q-\widetilde \Pi_{P_0} q| \le C \, h \,
\|q\|_2
\end{equation}
for all $\V \in \H^1$  and $q \in H^2$.

%

\noindent  We introduce  the finite element spaces
\begin{eqnarray}
\label{eq:defp1nc}
  P^d_1 &=& \{ v \in L^2 \; ; \; \forall \, K \in {\mathcal{T}}_h, \; \; v|_K \hbox{ is affine} \} \, , \nonumber\\
  P^{nc}_1 &=& \{ v_h \in P^d_1 \; ; \; \forall \, \sigma \in {\mathcal{E}}^{int}_h, \,
  v_h|_{K_\sigma}(\x_\sigma)=v_h|_{L_\sigma}(\x_\sigma) \, , \nonumber \\
  \P^{c}_1 &=& \{ \V_h \in (P^{d}_1)^2 \, ; \; \V_h \hbox{ is continuous} \hspace{.1cm} \hbox{ and } \hspace{.1cm}
  \V_h|_{\partial \Omega}=\mathbf{0}\}. \nonumber
 \end{eqnarray}
We have $\P^c_1 \subset \H^1_0$. We define the projection operator
$\Pi_{\P^{c}_1}:\H^1_0 \to \P^{c}_1$. For all $\V=(v_1,v_2) \in
\H^1_0$,
 $\Pi_{\P^{c}_1} \V=(v^1_h,v^2_h) \in \P^{c}_1$ is given by
\begin{equation*}
  \forall \, \hbox{\boldmath $\phi$ \unboldmath}_{\hspace{-1.5mm} h}=(\phi^1_h,\phi^2_h)  \in \P^c_1 \, ,\hspace{1cm}   \sum_{i=1}^2 \big( \nabla v^i_h,\nabla \phi^i_h)
  = \sum_{i=1}^2 \big( \nabla v_i,\nabla \phi^i_h).
\end{equation*}
The operator $\Pi_{\P^c_1}$
 is stable for the  $\H^1$ norm and
 (\cite{brenns} p. 110) there exists  $C>0$ such that for all $\V
\in \H^1$
\begin{equation}
\label{eq:errintp1c}
   |\V-\Pi_{\P^c_1} \V| \le C \, h \, \|\V\|_1.
\end{equation}
Let us address now the space  $P^{nc}_1$. If $q_h \in P^{nc}_1$, we
have usually $\nabla q_h \not \in \L^2$. Thus we define the operator
$\widetilde \nabla_h: P^{nc}_1 \hspace{-.2cm}\to \P_0 $
 by setting  for all $q_h \in P^{nc}_1$ and all   $K \in {\mathcal{T}}_h$
\begin{equation*}
\label{eq:defg}
  \widetilde \nabla_h q_h |_K = \frac{1}{|K|} \, \int_K \nabla q_h \, d\x.
\end{equation*}
The associated norm is given by
\begin{equation*}
\label{eq:defnormp1nc}
  \|q_h\|_{1,h}=\left(|q_h|^2+|\widetilde \nabla_h q_h|^2 \right)^{1/2}.
\end{equation*}
We also have a  Poincar\'e inequality: there exists  $C>0$ such that
for all $q_h \in P^{nc}_1 \cap L^2_0$
\begin{equation}
\label{eq:inpoinp1nc}
  |q_h| \le C \, |\widetilde \nabla_h q_h|.
\end{equation}
We define the projection operator $\Pi_{P^{nc}_1}$. For all $q_h \in
P^{nc}_1$, $\Pi_{P^{nc}_1} q_h$ is given by
\begin{equation}
\label{eq:defpp1nc}
   \forall \, \phi \in L^2 \, , \hspace{1cm} (\Pi_{P^{nc}_1} q_h,\phi)=(q_h,\phi).
\end{equation}
We have the following result.
\begin{prop}
\label{propp1ncu}

\noindent If $q_h \in P_0$,  $\Pi_{P^{nc}_1} q_h$ is given by
\begin{eqnarray*}
 \forall \, \sigma \in {\mathcal{E}}^{int}_h, \; \; \;(\Pi_{P^{nc}_1} q_h)(\x_\sigma) &=& \frac{|K_\sigma|}{|K_\sigma| + |L_\sigma|} \, q_{K_\sigma}
 +\frac{|L_\sigma|}{|K_\sigma| + |L_\sigma|} \, q_{L_\sigma}, \label{eq:p1ncta} \\
\forall \, \sigma \in {\mathcal{E}}^{ext}_h, \; \; \;(\Pi_{P^{nc}_1}
q_h)(\x_\sigma) &=&  q_{K_\sigma}. \label{eq:p1nctb}
\end{eqnarray*}
\end{prop}
\noindent {\sc Proof.} For all edge $\sigma \in {\mathcal{E}}_h$, we
define the function  $\psi_\sigma \in P^{nc}_1$ by setting
\begin{equation*}
  \psi_\sigma(\x_{\sigma'})= \left\{
    \begin{array}{l}
      1 \hbox{ if } \sigma=\sigma',  \\
      0 \hbox{ otherwise}.
    \end{array}
  \right.
.
\end{equation*}
Let us notice that $\psi_\sigma$ vanishes outside $K_\sigma \cup
L_\sigma$ if $\sigma \in {\mathcal{E}}^{int}_h$ and outside
$K_\sigma$ if $\sigma \in {\mathcal{E}}^{ext}_h$. 
Let   $\sigma \in {\mathcal{E}}^{int}_h$. Using a quadrature formula
we get
\begin{equation*}
\label{eq:projp1a} (\Pi_{P^{nc}_1} q_h, \psi_\sigma)=\left(
\frac{|K_\sigma|}{3} + \frac{|L_\sigma|} {3}\right)(\Pi_{P^{nc}_1}
q_h)(\x_\sigma)
\end{equation*}
and
\begin{equation*}
(q_h,\psi_\sigma)= q_{K_\sigma} \, \frac{|K_\sigma|}{3} +
q_{L_\sigma} \, \frac{|L_\sigma|}{3}.
\end{equation*}
 For an edge $\sigma \in {\mathcal{E}}^{ext}_h$ we have
$
\label{eq:projp1c} (\Pi_{P^{nc}_1} q_h, \psi_\sigma)
=\frac{|K_\sigma|}{3} \, (\Pi_{P^{nc}_1} q_h)(\x_\sigma)
$
and
$
\label{eq:projp1d}
  (q_h,\psi_\sigma)
  = q_{K_\sigma} \, \frac{|K_\sigma|}{3}
$.
 By plugging these equations into (\ref{eq:defpp1nc}) with $\phi=\psi_\sigma$, we get the
result. \qed

\noindent  We finally introduce the Raviart-Thomas spaces
\begin{align*}
\label{eq:defrt0}
  \mathbf{RT^d_0}=& \{ \V_h \in \P^d_1 \; ; \quad \forall \, \sigma \in {\mathcal{E}}_K,
   \quad \V_h|_K \cdot \N_{K,\sigma} \hbox{ is a  constant,}  \hspace{.1cm} \hbox{ and }  \hspace{.1cm}
    \V_h \cdot \N|_{\partial \Omega}=0  \} \, , \nonumber \\
  \mathbf{RT_0}=& \{ \V_h \in \mathbf{RT^d_0} \; ; \quad \forall \, K \in {\mathcal{T}}_h,
   \quad \forall \, \sigma \in {\mathcal{E}}_K,
 \quad   \V_h|_{K_\sigma} \cdot \N_{K_\sigma,\sigma} = \V_h|_{L_\sigma} \cdot \N_{K_\sigma,\sigma}\}.
\end{align*}
For all $\V_h \in \mathbf{RT_0}$, $K \in {\mathcal{T}}_h$ and
$\sigma \in {\mathcal{E}}_K$ we set $(\V_h \cdot
\N_{K,\sigma})_\sigma=\V_h|_K \cdot \N_{K,\sigma}$. \noindent We
define   the operator $\Pi_{\mathbf{RT_0}}:\H^1 \to \mathbf{RT_0}$.
For all $\V \in \H^1$, $\Pi_{\mathbf{RT_0}} \V \in \mathbf{RT_0}$ is
given by
\begin{equation}
\label{eq:defprt0} \forall \, K \in {\mathcal{T}}_h \, ,
\hspace{.5cm}
  \forall \, \sigma\in {\mathcal{E}}_K \, , \hspace{1cm}
    (\Pi_{\mathbf{RT_0}} \V \cdot \N_{K,\sigma})_\sigma=\frac{1}{|\sigma|}\int_\sigma \V \, d\sigma.
\end{equation}
One checks  \cite{brezzfor} that  there exists $C>0$ such that for all
$\V \in \H^1$
\begin{equation}
\label{eq:esterrirt0}
  |\V-\Pi_{\mathbf{RT_0}} \V| \le C \, h \, \|\V\|_1.
\end{equation}
The following result will be  useful.
\begin{prop}
\label{prop:rt0div}
  For all $\V \in \H^1$ such that $\hbox{\rm div} \, \V=0$, we have $\Pi_{\mathbf{RT_0}} \V \in \P_0$.
\end{prop}
\noindent {\sc Proof.} Let $\V_h=\Pi_{\mathbf{RT_0}} \V$ and $K \in
{\mathcal{T}}_h$.  According to
 \cite{brezzfor} there exists $\mathbf{a}_K \in \mathbb{R}^2$
and $b_K \in \mathbb{R}$ such that: $ \forall \, \x \in K \, ,
\hspace{.1cm} \V_h(\x)=\mathbf{a}_K+b_K \, \x$. Thus $\hbox{div} \,
\V_h|_K=2 \, b_K$. On the other hand, according to the divergence
formula and (\ref{eq:defprt0})
\begin{equation*}
\label{eq:rt0div1}
  0=\int_K \hbox{\rm div} \,\V \, d\x=\int_{\partial K} \V \cdot \N \, d\gamma
=\int_{\partial K} \V_h \cdot \N \, d\gamma=\int_K \hbox{\rm div}
\,\V_h \, d\x.
\end{equation*}
Hence $b_K=0$ and we get: $ \forall \, \x \in K \, , \hspace{.2cm}
\V_h(\x)=\mathbf{a}_K$.  \qed

\subsection{The discrete operators}
\label{subsec:opd}

The equations (\ref{eq:mom})--(\ref{eq:incomp}) use the differential
operators gradient, divergence and laplacian.
 Using the spaces of section \ref{subsec:espd} we define their discrete counterparts.
  The discrete gradient $ \nabla_h: P_0 \to \P_0$ is built using a
linear interpolation on the edges of the mesh (see \cite{zimm} for
details). This kind of construction has also be considered in \cite{egh2}. We set for all $q_h \in P_0$ and all  $K \in
{\mathcal{T}}_h$
\begin{eqnarray}
\label{eq:defg}
  \nabla_h \, q_h|_K&=&\frac{1}{|K|} \sum_{\sigma \in {\mathcal{E}}_K \cap {\mathcal{E}}^{int}_h}
|\sigma| \, \Big( \alpha_{K_\sigma,L_\sigma} \, q_{K_\sigma} +\alpha_{L_\sigma,K_\sigma} \, q_{L_\sigma} \Big) \, \N_{K,\sigma}  \nonumber \\
&+&\frac{1}{|K|}\sum_{\sigma \in {\mathcal{E}}_K \cap
{\mathcal{E}}^{ext}_h} |\sigma| \, q_{K_\sigma} \, \N_{K,\sigma}.
\end{eqnarray}
We have the following result \cite{zimm}.
\begin{prop}
\label{prop:kergradh}
  If $q_h \in L^2_0$ is such that $\nabla_h q_h=0$, then $q_h=0$.
\end{prop}
\noindent  The discrete divergence operator $\hbox{div}_h: \P_0 \to
P_0$ is built so that it is adjoint to the operator $\nabla_h$
(proposition \ref{prop:propadjh} below). We set for all $q_h \in
P_0$ and all  $K \in {\mathcal{T}}_h$
\begin{equation}
\label{eq:defd} \hbox{div}_h \, \V_h|_K=\frac{1}{|K|} \sum_{\sigma
\in {\mathcal{E}}_K \cap {\mathcal{E}}^{int}_h} |\sigma| \, \Big(
\alpha_{L_\sigma,K_\sigma} \, \V_{K_\sigma}
+\alpha_{K_\sigma,L_\sigma} \, \V_{L_\sigma} \Big)
 \cdot \N_{K,\sigma}.
\end{equation}
The first discrete laplacian $\Delta_h:P_0 \to P_0$ ensures that the
incompressibility constraint (\ref{eq:incomp}) is satisfied in a
discrete sense (proposition \ref{prop:umrt0}). We set for all $q_h
\in P_0$
\begin{equation}
\label{eq:deflap} \Delta_h q_h=\hbox{div}_h (\nabla_h q_h).
\end{equation}
\noindent The second discrete laplacian
$ \Deltat_h: \P_0 \to \P_0$ is the usual operator in finite volume
schemes \cite{eymgal}. We set for all $\V_h \in \P_0$ and all $K \in
{\mathcal{T}}_h$
\begin{equation*}
\label{eq:defl}
  \Deltat_h \V_h|_K
  = \frac{1}{|K|} \sum_{\sigma \in {\mathcal{E}}_K \cap {\mathcal{E}}^{int}_h}
 \tau_\sigma \, (\V_{L_\sigma} - \V_{K_\sigma} )
-\frac{1}{|K|}\sum_{\sigma \in {\mathcal{E}}_K \cap
{\mathcal{E}}^{ext}_h}
 \tau_\sigma \, \V_{K_\sigma}.
\end{equation*}


\vspace{.1cm}

\noindent In order to approximate the convection term  $(\U \cdot
\Nabla) \U$ in (\ref{eq:mom}) we define a bilinear form $\bt_h: \P_0
\times \P_0 \to \P_0$ using the well-known  upwind scheme (\cite{eymgal} p. 766). For all
$\U_h \in \P_0$, $\V_h \in \P_0$, and all  $K \in {\mathcal{T}}_h$
we have
\begin{equation}
\label{eq:defbth}
  \bt_h(\U_h,\V_h)\big|_K=
\frac{1}{|K|} \sum_{\sigma \in {\mathcal{E}}_K \cap
{\mathcal{E}}^{int}_h} \,
   |\sigma| \, \Big( (\U_\sigma \cdot \N_{K,\sigma})^+ \, \V_K +
(\U_\sigma \cdot \N_{K,\sigma})^- \, \V_{L_\sigma}  \Big).
\end{equation}
We have set $\U_\sigma=\alpha_{L_\sigma,K_\sigma} \, \U_{K_\sigma}
+\alpha_{K_\sigma,L_\sigma} \, \U_{L_\sigma}$ and $a^+=\max(a,0)$,
$a^-=\min(a,0)$ for all $a\in\mathbb{R}$.  Lastly, we define the
trilinear form $ \b_h: \P_0 \times \P_0 \times \P_0 \to
\mathbb{R}^2$ as follows. For all $\U_h \in \P_0$, $\V_h \in \P_0$,
$\W_h \in \P_0$, we set
\begin{equation}
\label{eq:defbbh}
  \b_h(\U_h,\V_h,\W_h)=\sum_{K \in {\mathcal{T}}_h} |K| \, \W_K \cdot
   \bt_h(\U_h,\V_h)\big|_K .
\end{equation}

\section{The scheme}
\label{sec:presschema}

We have defined in section \ref{sec:notd} the discretization in
space.  We now have  to define a discretization in time, and treat
the incompressibility constraint (\ref{eq:incomp}). We use a
projection method to this end. This kind of method has been
introduced by {\sc Chorin} \cite{chorin}
 and {\sc Temam} \cite{temam}. The basic idea is the following.
The time interval $[0,T]$ is split with a time step $k$:
$[0,T]=\bigcup_{n=0}^N [t_n,t_{n+1}]$ with $N \in \mathbb{N}^*$ and
$t_n=n \, k$ for all $n\in\{0,\dots,N\}$. For all
$m\in\{2,\dots,N\}$, we compute  (see equation (\ref{eq:mombdf})
below) a first velocity field $\Ut^m_h \simeq \U(t_m)$ using only
equation (\ref{eq:mom}). We use a second-order BDF scheme for the
discretization in time. We then project $\Ut^m_h$ (see equation
(\ref{eq:projib}) below) over a subspace of $\P_0$. We get a a
pressure field  $p^m_h \simeq p(t_m)$ and a second velocity field
$\U^m_h \simeq \U(t_m)$, which fulfills the incompressibilty
constraint (\ref{eq:incomp}) in a discrete sense. The algorithm goes
as follows.

\noindent First, for all $m\in\{0,\dots,N\}$, we set
$\F^m_h=\Pi_{\P_0} \F(t_m)$. Since the operator $\Pi_{\P_0}$ is
stable for the $\L^2$-norm we get
\begin{equation}
\label{eq:propfh}
  |\F^m_h|=|\Pi_{\P_0} \F(t_m)| \le |\F(t_m)| \le \|\F\|_{{\mathcal{C}}(0,T;\L^2)}.
\end{equation}
We start with the initial values
\begin{equation*}
  \U^0_h \in \P_0 \cap \mathbf{RT_0} \, , \hspace{1cm}   \U^1_h \in \P_0 \cap \mathbf{RT_0}
  \,\hspace{1cm}   p^1_h \in P_0 \cap L^2_0.
\end{equation*}
For all $n \in \{1,\dots,N\}$, $(\Ut^{n+1}_h,p^{n+1}_h,\U^{n+1}_h)$
is deduced from $(\Ut^n_h,p^n_h,\U^n_h)$ as follows.

\begin{itemize}
\item $\Ut^{n+1}_h \in \P_0$ is given by
\begin{equation}
\label{eq:mombdf} \frac{3 \, \Ut^{n+1}_h -4 \, \U^n_h+\U^{n-1}_h}{2
\, k} - \frac{1}{\hbox{Re}} \, \Deltat_h \Ut^{n+1}_h+\bt_h(2\,
\U^n_h-\U^{n-1}_h, \Ut^{n+1}_h)+\nabla_h p^n_h =\F^{n+1}_h  \, ,
\end{equation}

\item $p^{n+1}_h \in P^{nc}_1 \cap L^2_0$ is the solution of
\begin{equation}
\label{eq:projia} \Delta_h (p^{n+1}_h-p^n_h)=\frac{3}{2 \, k}\,
\hbox{div}_h \, \Ut^{n+1}_h  \, ,
\end{equation}

\item $\U^{n+1}_h \in  \P_0$ is deduced by
\begin{equation}
\label{eq:projib}
 \U^{n+1}_h = \Ut^{n+1}_h -\frac{2\, k}{3} \, \nabla_h (p^{n+1}_h - p^n_h).
\end{equation}
\end{itemize}
Existence and unicity of a solution to equation (\ref{eq:mombdf})
   is classical (\cite{eymgal} for example). Let us show
  that equation (\ref{eq:projia}) has also a unique solution.
Let  $q_h \in P_0 \cap L^2_0$ such that $\Delta_h q_h=0$.
 According to proposition \ref{prop:propadjh} we have for all $q_h \in P_0$
\begin{equation*}
  -(\Delta_h q_h,q_h)=-\big(\hbox{div}_h(\nabla_h q_h),q_h\big)=(\nabla_h q_h,\nabla_h q_h)=|\nabla_h q_h|^2.
\end{equation*}
Therefore we have  $\nabla_h q_h=0$. Using proposition
\ref{prop:kergradh} we get $q_h=0$. We have thus proved the unicity
of a solution for equation (\ref{eq:projia}). It is also the case
for the associated linear system.
 It implies that this linear system has indeed a solution. Hence it is also the case for
 equation  (\ref{eq:projia}).
\noindent Let us now prove  that for all $m\in\{0,\dots,N\}$,
$\U^m_h$ fulfills   (\ref{eq:incomp}) in a discrete sense.

\begin{lem}
\label{propdivrt0}

\noindent If $\V_h \in \mathbf{RT_0} \cap \P_0$ then
$ \hbox{\rm div}_h \,  \V_h=0$.
\end{lem}
\noindent {\sc Proof.} Let $K \in {\mathcal{T}}_h$. Since $\V_h \in
\mathbf{RT_0}$, definition (\ref{eq:defd}) reads
\begin{equation*}
  \hbox{div}_h \, \V_h|_K
=\frac{1}{|K|}\sum_{\sigma \in {\mathcal{E}}_K}
  |\sigma| \, (\alpha_{L_\sigma,K} +\alpha_{K,L_\sigma}) \, \V_K \cdot \N_{K,\sigma}.
\end{equation*}
Since $\alpha_{K_\sigma,L_\sigma}+\alpha_{L_\sigma,K_\sigma}=1$ we
conclude that
\begin{equation*}
\hspace{.9cm} \hbox{div}_h \, \V_h|_K= \frac{1}{|K|}\sum_{\sigma \in
{\mathcal{E}}_K}
  |\sigma| \, \V_K \cdot \N_{K,\sigma}
= \V_K \cdot \left( \frac{1}{|K|}\sum_{\sigma \in {\mathcal{E}}_K}
   |\sigma| \, \N_{K,\sigma} \right) =0. \hspace{.9cm} \qed
\end{equation*}

\begin{prop}
\label{prop:umrt0} For all $m\in\{0,\dots,N\}$ we have   $\hbox{\rm
div}_h \, \U^m_h=0$.
\end{prop}
\noindent {\sc Proof.} For $m\in\{0,1\}$ we have  $\U^0_h \in \P_0
\cap {\mathbf{RT_0}}$ and $\U^1_h \in \P_0 \cap {\mathbf{RT_0}}$.
Applying the lemma above we get the result.
 If $m\in\{2,\dots,N\}$, we apply the operator $\hbox{div}_h$ to
(\ref{eq:projia}) and compare with (\ref{eq:projib}).\qed

%


\section{Properties of the discrete operators}
\label{sec:propopd}

We prove that the differential operators in
(\ref{eq:mom})--(\ref{eq:incomp}) and the operators defined in
section \ref{subsec:opd} share similar properties.

\subsection{Properties of the discrete convective term}
\label{subsec:propconvh}

We define   $\bt: \H^1 \times \H^1 \to \L^2$. For all
 $\U \in \H^1$ and $\V=(v_1,v_2) \in \H^1$ we set
\begin{equation} \label{eq:defbt}
 \bt(\U,\V)= \big(\hbox{div}(v_1 \, \U),\hbox{div}(v_2 \, \U)\big).
\end{equation}
We show that the operator $\bt_h$ is a consistent approximation of
$\bt$.
\begin{prop}
\label{prop:consbh} There exists a constant $C>0$ such that for all
$\V \in \H^{2}$ and all $\U \in \H^2 \cap \H^1_0$ satisfying
$\hbox{\rm div} \, \U=0$
\begin{equation*}
\label{eq:consbh}
  \|\Pi_{\P_0}\bt(\U,\V) - \bt_h(\Pi_{\mathbf{RT_0}} \U, \widetilde \Pi_{\P_0} \V)\|_{-1,h}
\le C \, h  \,  \|\U\|_2 \,  \|\V\|_{1}.
\end{equation*}
\end{prop}
\noindent {\sc Proof.}  Let  $\U_h=\Pi_{\mathbf{RT_0}} \U$ and
$\V_h=\widetilde \Pi_{\P_0} \V$.
 According to  proposition \ref{prop:rt0div} we have $\U_h \in \P_0$.
%
Let $K \in {\mathcal{T}}_h$. According to the  divergence formula
and (\ref{eq:defprojp0}) we have
\begin{equation*}
\label{eq:eqconsb2}
 \Pi_{\P_0}\bt(\U,\V) |_K
 =\frac{1}{|K|} \sum_{\sigma \in {\mathcal{E}}_K \cap {\mathcal{E}}^{int}_h} \int_\sigma \V \, (\U \cdot \N) \, d\sigma.
\end{equation*}
On the other hand, let us rewrite $\bt_h(\U_h,\V_h)$. Let $\sigma
\in {\mathcal{E}}_K \cap {\mathcal{E}}^{int}_h$. Setting
\begin{equation*}
  \V_{K,L_\sigma}= \left\{
\begin{array}{ll}
   \V_{K} &\hbox{ si } (\U_h \cdot \N_{K,\sigma})_\sigma \ge 0 \\
   \V_{L_\sigma} &\hbox{ si } (\U_h \cdot \N_{K,\sigma})_\sigma < 0
\end{array}
\right.
\end{equation*}
one  checks that
  $\V_K \, (\U_\sigma \cdot \N_{K,\sigma})^+ + \V_{L_\sigma} \,  (\U_\sigma \cdot \N_{K,\sigma})^-
=  \V_{K,L_\sigma} \, (\U_\sigma \cdot \N_{K,\sigma})$.
By definition
  $\U_\sigma \cdot \N_{K,\sigma}
   =\alpha_{L_\sigma,K} \, ( \U_K \cdot \N_{K,\sigma})+
   \alpha_{K,L_\sigma} \, (\U_{L_\sigma} \cdot \N_{K,\sigma})$ ;
since $\U_h \in \mathbf{RT_0}$ we get
   $\U_\sigma \cdot \N_{K,\sigma}
   = (\alpha_{L_\sigma,K}+\alpha_{K,L_\sigma}) \, (\U_K \cdot \N_{K,\sigma})
   =(\U_K \cdot \N_{K,\sigma})=(\U_h \cdot \N_{K,\sigma})_\sigma$.
Using at last (\ref{eq:defprt0}),  we deduce from (\ref{eq:defbth})
\begin{equation*}
\label{eq:eqconsb3} \bt_h(\U_h,\V_h)|_K=\frac{1}{|K|} \sum_{\sigma
\in {\mathcal{E}}_K \cap {\mathcal{E}}^{int}_h} \int_\sigma
\V_{K,L_\sigma} \, (\U \cdot \N_{K,\sigma}) \, d\sigma.
\end{equation*}
Thus
\begin{equation*}
\label{eq:eqconsb4}
  \big( \Pi_{\P_0}\bt(\U,\V)-\bt_h(\U_h,\V_h) \big)|_K
 =\frac{1}{|K|} \sum_{\sigma \in {\mathcal{E}}_K \cap {\mathcal{E}}^{int}_h}
 \int_\sigma (\V-\V_{K,L_\sigma}) \, (\U \cdot \N) \, d\sigma.
\end{equation*}
 Let $\psig_h \in \P_0$. We have
\begin{eqnarray}
\label{eq:eqconsb5}
 \hspace{-.1cm} \Big( \Pi_{\P_0}\bt(\U,\V) - \bt_h(\U_h, \V_h), \psig_h \Big)
&=&\sum_{K \in {\mathcal{T}}_h}   \psig_K  \sum_{\sigma \in
{\mathcal{E}}_K \cap {\mathcal{E}}^{int}_h}
 \int_\sigma (\V-\V_{K,L_\sigma}) \, (\U \cdot \N) \, d\sigma \nonumber \\
 &=& \sum_{\sigma \in {\mathcal{E}}^{int}_h} (\psig_{K_\sigma}-\psig_{L_\sigma})  \int_\sigma (\V-\V_{K_\sigma,L_\sigma}) \, (\U \cdot \N) \, d\sigma.
\end{eqnarray}
Let $\sigma \in {\mathcal{E}}^{int}_h$. We want to estimate the
integral over $\sigma$. Since we work in a two-dimensional domain,
we have the Sobolev injection $\H^{2} \subset \L^\infty$.
 Thus
\begin{equation*}
\label{eq:eqconsb6}
   \left| \int_\sigma (\V-\V_{K_\sigma,L_\sigma}) \, (\U \cdot \N) \, d\sigma \right|
   \le \|\U\|_{\L^\infty}  \int_\sigma |\V-\V_{K_\sigma,L_\sigma}| \, d\sigma
   \le C \, \|\U\|_{2}  \int_\sigma |\V-\V_{K_\sigma,L_\sigma}| \, d\sigma.
\end{equation*}
Let us first assume that  $\V\in {\mathcal{C}}^1$. We set
\begin{equation*}
  \x_{K_\sigma,L_\sigma}= \left\{
\begin{array}{ll}
   \x_{K_\sigma} &\hbox{ si } (\U_h \cdot \N_{K,\sigma})_\sigma \ge 0 \\
   \x_{L_\sigma} &\hbox{ si } (\U_h \cdot \N_{K,\sigma})_\sigma < 0
\end{array}
\right. .
\end{equation*}
 If $\x \in \sigma$, we have the following Taylor expansion
\begin{equation*}
\V(\x)-\V_{K_\sigma,L_\sigma}=\V(\x)-\V(\x_{K_\sigma,L_\sigma})=\int^1_0
\Nabla \V \, (t \, \x+(1-t) \, \x_{K_\sigma,L_\sigma}) \,
(\x-\x_{K_\sigma,L_\sigma}) \, dt.
\end{equation*}
We have $|\x-\x_{K_\sigma,L_\sigma}| \le h$. Thus, integrating over
$\sigma$ and using the Cauchy-Schwarz inequality, we get
\begin{equation*}
\int_\sigma |\V-\V_{K_\sigma,L_\sigma}| \, d\sigma \le \sqrt{2}
\left(  \int_\sigma \int^1_0 |\Nabla \V \, (t \, \x+(1-t) \,
\x_{K_\sigma,L_\sigma})|^2  \, h \, \sqrt{t} \, dt \, d\sigma
\right)^{1/2}.
\end{equation*}
We then use the change of variable $(t,\x) \to \mathbf{y}=t \,
\x+(1-t) \, \x_{K_\sigma,L_\sigma}$.
 Let $D_\sigma$ be the quadrilateral domain given by the endpoints of  $\sigma$, $\x_{K_\sigma}$ and $\x_{L_\sigma}$.
   The domain $[0,1] \times \sigma$ becomes $D_{K_\sigma,L_\sigma}$
   with
\begin{equation*}
  D_{K_\sigma,L_\sigma}= \left\{
\begin{array}{rl}
   D_\sigma \cap {K_\sigma} &\hbox{ si } (\U_h \cdot \N_{K,\sigma})_\sigma \ge 0 \\
   D_\sigma \cap {L_\sigma} &\hbox{ si } (\U_h \cdot \N_{K,\sigma})_\sigma < 0
\end{array}
\right. .
\end{equation*}
For all $t \in [0,1]$ we have  $h\sqrt{t} \le h \, t \le C \,
d(\x_{K_\sigma,L_\sigma},\sigma) \, t$ thanks to the hypothesis  on
the mesh. We check easily that $d(\x_{K_\sigma,L_\sigma},\sigma) \,
t \, dt \, d\sigma=d\mathbf{y}$.
 Thus we get
\begin{equation*}
\label{eq:eqconsb7c} \int_\sigma |\V-\V_{K_\sigma,L_\sigma}| \,
d\sigma \le
 C \, h\left(  \int_{D_{K_\sigma,L_\sigma}} |\Nabla \V \, (\mathbf{y})|^2  \,  d\mathbf{y} \right)^{1/2}.
\end{equation*}
Since $({\mathcal{C}}^1)^2$ is dense in $\H^2$, this estimate  still
holds for  $\V \in \H^2$. Plugging this estimate into
(\ref{eq:eqconsb5}) and using the Cauchy-Schwarz inequality we get
\begin{eqnarray*}
&& \left|  \big( \Pi_{\P_0}\bt(\U,\V) - \bt_h(\Pi_{\mathbf{RT_0}}
\U, \widetilde \Pi_{\P_0} \V), \psig_h \big) \right|
  \\
&& \le C \, h \, \|\U\|_{\H^{2}} \left( \sum_{\sigma \in
{\mathcal{E}}^{int}_h} |\psig_{L_\sigma} - \psig_{K_\sigma}|^2
\right)^{1/2}
 \left( \sum_{\sigma \in {\mathcal{E}}^{int}_h}   \int_{D_{K_\sigma,L_\sigma}} |\Nabla \V \, (\mathbf{y})|^2
  \,  d\mathbf{y} \right)^{1/2}
\end{eqnarray*}
so that
$\left|  \big( \Pi_{\P_0}\bt(\U,\V) - \bt_h(\Pi_{\mathbf{RT_0}} \U,
\widetilde \Pi_{\P_0} \V), \psig_h \big) \right|
  \le C \, h \,  \|\U\|_{\H^{2}} \, \|\psig_h\|_{1,h} \, \|\V\|_1$.
Using then definition (\ref{eq:defnormdp0}), we get
the result. \qed

\noindent Let us consider now the operator $\b_h$. Let $\U \in \H^1$ and
$\V \in \L^\infty \cap \H^1$ with $\hbox{div} \, \U \ge 0$.
 Integrating by parts we deduce from (\ref{eq:defbt}):
    $\int_\Omega \V \cdot \bt(\U,\V) \, d\x=\int_\Omega  \frac{|\V|^2}{2}\,  \hbox{div} \, \U\, d\x \ge 0$.
The discrete operator $\b_h$ shares a similar property.

\begin{prop}
\label{prop:posbh} Let $\U_h \in \P_0$ such that $\hbox{\rm div}_h
\, \U_h \ge  0$. For all $\V_h \in \P_0$  we have
\begin{equation*}
\label{eq:bhpos}
  \b_h(\U_h,\V_h,\V_h) \ge 0.
\end{equation*}
\end{prop}
\noindent {\sc Proof.} Remember that for all edges $\sigma \in
{\mathcal{E}}^{int}_h$, two triangles $K_\sigma$ et $L_\sigma$ share
$\sigma$ as an edge. We denote by $K_\sigma$ the one such that
$\U_{\sigma} \cdot \N_{K_\sigma,\sigma} \ge 0$. Using the algebraic
identity $2 \, a \, (a-b)=a^2-b^2+(a-b)^2$ we deduce from
(\ref{eq:defbbh})
\begin{eqnarray*}
  2 \, \b_h(\U_h,\V_h,\V_h)&=&
2\sum_{\sigma \in {\mathcal{E}}^{int}_h} |\sigma| \, \V_{K\sigma}
\cdot (\V_{K_\sigma} - \V_{L_\sigma})
 \, (\U_{\sigma} \cdot \N_{K_\sigma,\sigma}) \\
&=&\sum_{\sigma \in {\mathcal{E}}^{int}_h} |\sigma| \,  \Big(
|\V_{K\sigma}|^2 - |\V_{L_\sigma}|^2+|\V_{K_\sigma} -
\V_{L_\sigma}|^2 \Big) \, (\U_{\sigma} \cdot \N_{K_\sigma,\sigma})
\end{eqnarray*}
so that
$2 \, \b_h(\U_h,\V_h,\V_h) \ge  \sum_{\sigma \in
{\mathcal{E}}^{int}_h} |\sigma|  \, \Big( |\V_{K\sigma}|^2 -
|\V_{L_\sigma}|^2 \Big) \, (\U_{\sigma} \cdot
\N_{K_\sigma,\sigma})$.
This sum can be written  as a sum over the triangles of the mesh. We
get
\begin{equation*}
2 \, \b_h(\U_h,\V_h,\V_h) \ge  \sum_{K \in {\mathcal{T}}_h}
|\V_{K_\sigma}|^2 \sum_{\sigma \in {\mathcal{E}}_K \cap
{\mathcal{E}}^{int}_h} |\sigma| \, (\U_\sigma \cdot
\N_{K_\sigma,\sigma}).
\end{equation*}
Using finally definition (\ref{eq:defd}) we get
\begin{equation*}
\hspace{2.3cm} 2 \, \b_h(\U_h,\V_h,\V_h) \ge \sum_{K \in
{\mathcal{T}}_h} |K| \, |\V_{K}|^2 \,    (\hbox{\rm div}_h \, \U_h)
|_K \ge 0. \hspace{2.3cm} \qed
\end{equation*}

\noindent The following result states that the operator $\b_h$ is
stable for suitable norms.

\begin{prop}
\label{prop:stabbth} There exists a constant $C>0$ such that for all
 $\V_h \in \P_0$, $\W_h\in \P_0$, $\U_h \in \P_0$ satisfying  $\hbox{ \rm div}_h \, \U_h=0$
\begin{equation*}
 |\b_h(\U_h,\V_h,\V_h)| \le C \, |\U_h| \,
\|\V_h\|_{h} \, \|\V_h\|_h.
\end{equation*}
\end{prop}
\noindent {\sc Proof.} For all triangle $K \in {\mathcal{T}}_h$ and
all edge  $\sigma \in {\mathcal{E}}_K \cap {\mathcal{E}}^{int}_h$,
we have
\begin{equation*}
  (\U_\sigma \cdot \N_{K,\sigma})^+ \, \V_K +  (\U_\sigma \cdot \N_{K,\sigma})^-
  \, \V_{L_\sigma}
  =(\U_\sigma \cdot \N_{K,\sigma}) \, \V_K
 - |(\U_\sigma \cdot \N_{K,\sigma})| \, (\V_{L_\sigma}-\V_K).
\end{equation*}
This way, we deduce from  (\ref{eq:defbh})
 $\b_h(\U_h,\V_h,\W_h) = S_1 +S_2$
with
\begin{eqnarray*}
 S_1 &=&
\sum_{K \in {\mathcal{T}}_h} \V_K \cdot \W_K \sum_{\sigma \in
{\mathcal{E}}_K \cap {\mathcal{E}}^{int}_h} |\sigma|
\,  (\U_\sigma \cdot \N_{K,\sigma}) \, , \nonumber \\
S_2 &=& -\sum_{K \in {\mathcal{T}}_h} \W_K \cdot \sum_{\sigma \in
{\mathcal{E}}_K \cap {\mathcal{E}}^{int}_h} |\sigma|\, |\U_\sigma
\cdot \N_{K,\sigma}| \, (\V_{L_\sigma}-\V_K) \label{eq:defs2}.
\end{eqnarray*}
By writing the sum over the edges as a sum over the triangles we get
\begin{equation*}
S_2 = -\sum_{\sigma \in  {\mathcal{E}}^{int}_h}   |\sigma|\,
|\U_\sigma \cdot \N_{K,\sigma}| \, (\V_{L_\sigma}-\V_K) \cdot
(\W_{L_\sigma}-\W_K).
\end{equation*}
Using the Cauchy-Schwarz inequality we get
\begin{equation*}
  |S_2| \le  h \, \|\U_h\|_{\infty} \, \left( \sum_{\sigma \in {\mathcal{E}}^{int}_h}
  |\V_{L_\sigma} - \V_{K_\sigma}|^2 \right)^{1/2} \,  \left( \sum_{\sigma \in {\mathcal{E}}^{int}_h} |\W_{L_\sigma} - \W_{K_\sigma}|^2 \right)^{1/2}.
\end{equation*}
Since $\U_h \in \P_0$ we have the inverse inequality \cite{eymgal}
  $h \, \|\U_h\|_{\infty} \le C \, |\U_h|$.
 Using (\ref{eq:proptaus}) and  (\ref{eq:defh1d}) we have
\begin{equation*}
 \sum_{\sigma \in {\mathcal{E}}^{int}_h} |\V_{L_\sigma} - \V_{K_\sigma}|^2
 \le  C \, \sum_{\sigma \in {\mathcal{E}}^{int}_h} \tau_\sigma \, |\V_{L_\sigma} - \V_{K_\sigma}|^2
 \le C \, \|\V_h\|^2_h
\end{equation*}
and
 $\sum_{\sigma \in {\mathcal{E}}^{int}_h} |\W_{L_\sigma} - \W_{K_\sigma}|^2  \le C \, \|\W_h\|^2_h$.
Therefore
  $|S_2| \le C \, |\U_h| \, \|\V_h\|_h \, \|\W_h\|_h$.
On the other hand we deduce from definition (\ref{eq:defd})
\begin{equation*}
\label{eq:estbh1} S_1 = \sum_{K \in {\mathcal{T}}_h}  |K| \,  (\V_K
\cdot \W_K) \, (\hbox{div}_h \, \U_h)|_K=0.
\end{equation*}
By combining the estimates for $S_1$ and $S_2$ we get the
result.\qed

\subsection{Properties of the discrete gradient}
\label{subsec:propgradd}
\begin{prop}
\label{prop:ininvgrad}There exists a constant $C>0$  such that for
all $q_h \in P_0$: $h \, |\nabla_h q_h| \le C \, |q_h|$.
\end{prop}
\noindent {\sc Proof.} Using (\ref{eq:defg}) and the Minkowski
inequality, we have for all triangle $K \in {\mathcal{T}}_h$
\begin{equation*}
  |K| \left| \nabla_h q_h \, |_K \right|^2
  \le \sum_{\sigma \in {\mathcal{E}}_K \cap {\mathcal{E}}^{int}_h}
  \frac{6 \, |\sigma|^2}{|K|} \, (q^2_{K}+q^2_{L_\sigma})
+ \sum_{\sigma \in {\mathcal{E}}_K \cap {\mathcal{E}}^{ext}_h}
  \frac{6 \, |\sigma|^2}{|K|} \,   q^2_{K}.
\end{equation*}
Let us sum over $K \in {\mathcal{T}}_h$. Since $|\sigma| \le h$,
using (\ref{eq:propairet}), we get
\begin{equation*}
  |\nabla_h q_h|^2 \le
  \frac{C}{h^2} \left( \sum_{K \in {\mathcal{T}}_h}
\sum_{\sigma \in {\mathcal{E}}_K \cap {\mathcal{E}}^{int}_h}
    \left(|K| \, q^2_{K}+|L_\sigma| \, q^2_{L_\sigma}\right)
+\sum_{K \in {\mathcal{T}}_h} \sum_{\sigma \in {\mathcal{E}}_K \cap
{\mathcal{E}}^{ext}_h} |K| \, q^2_{K} \right).
\end{equation*}
Thus
$ h^2 \, |\nabla_h q_h|^2 \le C  \sum_{K \in {\mathcal{T}}_h} |K| \,
q^2_{K} \le C \, |q_h|^2$. \qed
\vspace{.1cm}

\noindent We now prove that $\nabla_h$ is a consistent approximation
of the gradient.
\begin{prop}
\label{prop:consg}
There exists a constant $C>0$ such that for all $q \in H^2$
\begin{equation*}
     |\Pi_{\P_0}(\nabla q) -\nabla_h (\widetilde \Pi_{P_0} q)| \le C \, h \, \|q\|_2.
\end{equation*}
\end{prop}
\noindent {\sc Proof.} Let $K \in {\mathcal{T}}_h$. Using the
gradient formula and definition (\ref{eq:defg}) we get
\begin{equation*}
\label{eq:consg1}
  |K|  \left. \left( \Pi_{\P_0} (\nabla q) - \nabla_h (\widetilde \Pi_{P_0} q) \right)\right|_K
 = \int_K \nabla q \, d\x - |K| \left. \nabla_h (\widetilde \Pi_{P_0} q)  \right|_K
 = \sum_{\sigma \in {\mathcal{E}}_K} I^\sigma_{K}
\end{equation*}
where we have set for all edge  $\sigma \in {\mathcal{E}}_K \cap
{\mathcal{E}}^{int}_h$
\begin{equation*}
\label{eq:consg2} I^\sigma_K  =   \int_\sigma \Big( q -
\big(\alpha_{K,L_\sigma} \, q(\x_{K}) +\alpha_{L_\sigma,K}
 \, q(\x_{L_\sigma})\big)\Big) \, \N_{K,\sigma} \, d\sigma
\end{equation*}
and for all edge  $\sigma \in {\mathcal{E}}_K \cap
{\mathcal{E}}^{ext}_h$:
 $I^\sigma_K = \int_\sigma \big(q - q(\x_{K})\big)
\, \N_{K,\sigma} \, d\sigma$.
Squaring and using (\ref{eq:propairet}) we get
\begin{equation*}
|K|  \left| \left. \left( \Pi_{\P_0} (\nabla q) - \nabla_h
(\widetilde \Pi_{P_0} q) \right)\right|_K \right|^2 \le
\frac{3}{|K|} \sum_{\sigma \in {\mathcal{E}}_K} |I^\sigma_K|^2 \le
\frac{C}{h^2} \sum_{\sigma \in {\mathcal{E}}_K} |I^\sigma_K|^2.
\end{equation*}
Summing over the triangles $K \in {\mathcal{T}}_h$ we get
\begin{equation}
\label{eq:consg35}
   \left|  \Pi_{\P_0} (\nabla q) - \nabla_h (\widetilde \Pi_{P_0}   q)\right|^2
 \le  \frac{C}{h^2} \sum_{K \in {\mathcal{T}}_h} \sum_{\sigma
\in {\mathcal{E}}_K} |I^\sigma_K|^2 .
\end{equation}
We must estimate the  integral terms $I^\sigma_K$. Let $K \in
{\mathcal{T}}_h$. Let us first assume  that $q \in
{\mathcal{C}}^2(\overline \Omega)$. Let $\sigma \in {\mathcal{E}}_K
\cap {\mathcal{E}}^{int}_h$. For $\x \in \sigma$ we have the
following Taylor expansions
\begin{equation*}
\label{eq:consg4}
  q(\x_K)=q(\x) +\nabla q(\x) \cdot (\x_K-\x) + \int^1_0 \H(q) \, (t\x_K+(1-t)\x)
(\x_K-\x) \cdot (\x_K-\x)
  \, t \, dt \, ,
\end{equation*}
\begin{equation*}
\label{eq:consg5}
  q(\x_{L_\sigma})=q(\x) +\nabla q(\x) \cdot (\x_{L_\sigma}-\x) + \int^1_0 \H(q) \, (t\x_{L_\sigma}+(1-t)\x)
(\x_{L_\sigma}-\x) \cdot (\x_{L_\sigma}-\x)
  \, t \, dt \, ,
\end{equation*}
\begin{equation*}
\label{eq:consg6}
  \nabla q(\x)=\nabla q(\x_K)- \int^1_0  \Nabla \nabla q \, \big(t\x_K+(1-t)\x\big) (\x_K-\x) \; dt.
\end{equation*}
Plugging the last expansion into the two others  and integrating
over $\sigma$ we get
\begin{equation}
 \label{eq:consg65}
  \int_\sigma \big(q(\x_{K})-q\big) \, d\sigma= |\sigma| \, \nabla q(\x_K) \cdot (\x_K-\x_\sigma)
-A^\sigma_K+B^\sigma_K \, ,
\end{equation}
\begin{equation}
 \label{eq:consg7}
 \int_\sigma \big(q(\x_{L_\sigma})-q\big) \, d\sigma= |\sigma| \, \nabla q(\x_K) \cdot (\x_{L_\sigma}-\x_\sigma)
  -A^\sigma_{L_\sigma}+B^\sigma_{L_\sigma}.
\end{equation}
We have set for $T \in \{K_\sigma,L_\sigma\}$
\begin{equation}
\label{eq:defah}
 A^\sigma_T=\int_\sigma \int^1_0 \Nabla \nabla q \, (t\x_T+(1-t)\x)
 \, (\x_T-\x)  \, dt \, d\sigma \, ,
\end{equation}
\begin{equation}
\label{eq:defbh} B^\sigma_T=\int_\sigma \int^1_0 \H(q) \,
(t\x_T+(1-t)\x) (\x_T-\x) \cdot (\x_T-\x)  \, t \, dt \, d\sigma.
\end{equation}
One can bound these terms as  in the proof of proposition
\ref{prop:consbh}. We get
\begin{equation}
\label{eq:esttermcons}
 |A^\sigma_T|^2 \le C  \, h^2  \int_{D_\sigma} |\Nabla \nabla q \, (\mathbf{y})|^2  \,
 d\mathbf{y} \, ,
 \hspace{1cm} |B^\sigma_T|^2 \le C  \, h^4  \int_{D_\sigma} |\mathbf{H}(q)(\mathbf{y})|^2  \,
 d\mathbf{y}.
\end{equation}
 Now, let us multiply (\ref{eq:consg65}) by
$-\alpha_{K,L_\sigma} \, \N_{K,\sigma}$, (\ref{eq:consg7}) by
$-\alpha_{L_\sigma,K} \, \N_{K,\sigma}$ and sum  the equalities.
Since $\alpha_{L_\sigma,K}+\alpha_{K,L_\sigma}=1$ we have
\begin{eqnarray*}
&&   -\alpha_{L_\sigma,K}  \int_\sigma \big(q(\x_{K})-q\big) \,
\N_{K,\sigma} \, d\sigma
-\alpha_{K,L_\sigma} \int_\sigma \big(q(\x_{L_\sigma})-q\big) \, \N_{K,\sigma}\, d\sigma  \\
&& =
  \int_\sigma \Big( q -  \big(\alpha_{K_\sigma,L_\sigma} \, q(\x_{K,\sigma})
   +\alpha_{L_\sigma,K_\sigma} \, q(\x_{L,\sigma})\big)\Big) \,
   \N_{K,\sigma} \, d\sigma
=I^\sigma_K.
\end{eqnarray*}
On the other hand
\begin{equation*}
-\alpha_{K,L_\sigma} \, (\x_K-\x_\sigma) \cdot
\N_{K,\sigma}-\alpha_{L_\sigma,K} \,  (\x_{L_\sigma}-\x_\sigma)
\cdot \N_{K,\sigma} =-\alpha_{K,L_\sigma} \, \alpha_{L_\sigma,K} \,
(d_\sigma-d_\sigma)=0.
\end{equation*}
Therefore we get
$I^\sigma_K=-\alpha_{L_\sigma,K} \, \big(A^\sigma_K+B^\sigma_K\big)
\, \N_{K,\sigma}
 - \alpha_{K,L_\sigma} \, \big(A^\sigma_{L_\sigma}+B^\sigma_{L_\sigma}\big) \, \N_{K,\sigma}$.
Using estimates (\ref{eq:esttermcons}) we obtain
\begin{equation*}
\label{eq:consg105}
 |I^\sigma_K|^2 \le C  \, h^4  \int_{D_\sigma} (|\H(q)(\mathbf{y})|^2 + |\Nabla \nabla q(\mathbf{y})|^2)
 \, d\mathbf{y}.
\end{equation*}
We now consider  the case $\sigma \in {\mathcal{E}}_K \cap
{\mathcal{E}}^{ext}_h$. For $\x \in \sigma$ we have
\begin{equation*}
\label{eq:consg11}
  q(\x_K)=q(\x)\!+\!
  \nabla q(\x) \cdot (\x_{K}-\x)
 \!+\! \int^1_0 \! \! \H(q)  (t\x_K+(1-t)\x)
(\x_{K}-\x) \cdot (\x_{K}-\x)  t  dt.
\end{equation*}
Multiplying by $\N_{K,\sigma}$ and integrating over $\sigma$, we get
  $-I^\sigma_K= J^\sigma_K \, \N_{K,\sigma}+B^\sigma_K \,
  \N_{K,\sigma}$
with
$J^\sigma_K= \int_\sigma   \nabla q(\x) \cdot (\x_K-\x)  \, d\x$.
Since  $|\x_K-\x| \le h$ if $\x \in \sigma$, using a trace theorem,
we have
\begin{equation*}
\label{eq:consg12} |J^\sigma_K| \le C \, h^2 \, \|\nabla q
\|_{\L^\infty(\sigma)} \le C \, h^2 \left( \int_{D_\sigma}
\big(|\nabla q(\mathbf{y})|^2 +|\Nabla \nabla q(\mathbf{y})|^2\big)
\, d\mathbf{y} \right)^{1/2}.
\end{equation*}
By combining this estimate with (\ref{eq:esttermcons}), we get
\begin{eqnarray*}
\label{eq:consg13}
 |I^\sigma_K|^2 \le 2 \, |J^\sigma|^2+2\,|B^\sigma_K|^2 &\le&
C \, h^4 \,  \int_{D_\sigma} |\H(q)(\mathbf{y})|^2 \, d\mathbf{y} \nonumber \\
&+& C \, h^4 \, \int_{D_\sigma} (|\nabla q(\mathbf{y})|^2 +|\Nabla
\nabla q(\mathbf{y})|^2) \, d\mathbf{y}.
\end{eqnarray*}
The space ${\mathcal{C}}^2(\overline \Omega)$ is dense in $H^2$.
Therefore the bounds for  $I^\sigma_K$ still hold for $q \in H^2$.
Plugging these bounds into (\ref{eq:consg35}) we get the result.
\qed

\subsection{Properties of the discrete divergence}
\label{subsec:propdivh} The operators divergence and gradient are
adjoint:
 if $q \in H^1$ and $\V\in \H^1$ with $\V\cdot \N|_{\partial
 \Omega}=0$, we get $(\V,\nabla q)=-(q,\hbox{div} \, \V)$ by
 integrating by parts. For $\nabla_h$ and $\hbox{div}_h$ we state
\begin{prop}
\label{prop:propadjh}

\noindent If $\V_h \in \P_0$ and $q_h \in P_0$ we have: $
\label{eq:opdadj}
  (\V_h, \nabla_h q_h)=-(q_h,\hbox{\rm div}_h \, \V_h)
$.
\end{prop}
\noindent {\sc Proof.} Using   (\ref{eq:defg}) one checks that
 $(\V_h,\nabla_h q_h)=\sum_{K \in {\mathcal{T}}_h} q_K \,
 (S_1+S_2+S_3)$
with
\begin{equation*}
  S_1 = \sum_{\sigma \in {\mathcal{E}}_K \cap {\mathcal{E}}^{int}_h} |\sigma| \,  \alpha_{K,L_\sigma} \, \V_K \cdot \N_{K,\sigma} \, ,
  \hspace{1cm} S_2 = \sum_{\sigma \in {\mathcal{E}}_K \cap {\mathcal{E}}^{int}_h} |\sigma| \, \alpha_{K,L_\sigma} \, \V_{L_\sigma} \cdot
  \N_{L_\sigma,\sigma} \, ,
\end{equation*}
and
  $S_3 = \sum_{\sigma \in {\mathcal{E}}_K \cap {\mathcal{E}}^{int}_h} |\sigma|  \, \V_K \cdot \N_{K,\sigma}$.
Since $\alpha_{K,L_\sigma}+\alpha_{L_\sigma,K}=1$ we have
\begin{eqnarray*}
S_1&=&\sum_{\sigma \in {\mathcal{E}}_K \cap {\mathcal{E}}^{int}_h} |\sigma| \, (1-\alpha_{L_\sigma,K}) \,  \V_K \cdot \N_{K,\sigma} \\
&=&\sum_{\sigma \in {\mathcal{E}}_K \cap {\mathcal{E}}^{int}_h}
|\sigma| \, \V_K \cdot \N_{K,\sigma} -\sum_{\sigma \in
{\mathcal{E}}_K \cap {\mathcal{E}}^{int}_h} |\sigma| \,
\alpha_{L_\sigma,K} \, \V_K \cdot \N_{K,\sigma}.
\end{eqnarray*}
Since $\N_{L_\sigma,\sigma}=-\N_{K,\sigma}$, we also have
\begin{equation*}
S_2=\sum_{\sigma \in {\mathcal{E}}_K \cap {\mathcal{E}}^{int}_h}
|\sigma| \, \alpha_{K,L_\sigma} \, \V_{L_\sigma} \cdot
\N_{L_\sigma,\sigma} =-\sum_{\sigma \in {\mathcal{E}}_K \cap
{\mathcal{E}}^{int}_h} |\sigma| \, \alpha_{K,L_\sigma} \,
\V_{L_\sigma} \cdot \N_{K,\sigma}.
\end{equation*}
Therefore
\begin{equation*}
(\V_h,\nabla_h q_h) = -\sum_{\sigma \in {\mathcal{E}}_K \cap
{\mathcal{E}}^{int}_h} |\sigma| \, (\alpha_{L,K_\sigma} \, \V_K+
\alpha_{K,L_\sigma} \, \V_{L_\sigma}) \cdot \N_{K,L_\sigma}
+\sum_{\sigma \in {\mathcal{E}}_K} |\sigma|  \, \V_K \cdot
\N_{K,L_\sigma}.
\end{equation*}
Using definition (\ref{eq:defd}) we get
\begin{equation*}
(\V_h,\nabla_h q_h)  = - \sum_{K \in {\mathcal{T}}_h}|K| \, \hbox{div}_h \, \V_h|_K
+\sum_{K \in {\mathcal{T}}_h}\sum_{\sigma \in {\mathcal{E}}_K} |\sigma|  \, \V_K \cdot
\N_{K,L_\sigma}.
\end{equation*}
Since
$\sum_{\sigma \in {\mathcal{E}}_K} |\sigma| \, \V_K \cdot
\N_{K,L_\sigma}=\V_K \cdot \sum_{\sigma \in {\mathcal{E}}_K}
 |\sigma|  \, \N_{K,L_\sigma}=0$
we obtain finally
\begin{equation*}
  \hspace{1.9cm}(\V_h,\nabla_h q_h)= - \sum_{K \in {\mathcal{T}}_h} q_K \, |K| \, \hbox{div}_h \, \V_h|_K=-(q_h,\hbox{div}_h \, \V_h).
  \hspace{1.9cm} \qed
\end{equation*}
The divergence operator and the spaces $L^2_0$, $\H^1_0$ satisfy
 the following property, called inf-sup (or Babu\v{s}ka-Brezzi) condition (see \cite{giraultr} for example).
There exists a constant $C>0$ such that
\begin{equation}
\label{eq:infsc}
  \inf_{\begin{array}{l}
  {\scriptstyle q \in L^2_0 \backslash \{0\} }
  \end{array}}
\hspace{.2cm} \sup_{\V \in \H^1_0 \backslash \{\mathbf{0}\}}
-\frac{(q,\hbox{\rm div} \, \V)}{\|\V\|_1 |q|} \ge C.
\end{equation}
We will now prove that the operator $\hbox{\rm div}_h$ and the
spaces $P_0 \cap L^2_0$, $\P_0$ satisfy an analogous property. The
proof is based on the following lemma.
\begin{lem}
\label{lemma:gradmu}
  We assume that the mesh is uniform (i.e. the triangles of the mesh are equilateral). Then we have for all $q_h \in P_0$
  \begin{equation*}
  \label{eq:gradmu}
    \nabla_h q_h = \widetilde \nabla_h (\Pi_{P^{nc}_1} q_h).
  \end{equation*}
\end{lem}
\noindent {\sc Proof.} Since the mesh is uniform we have: $ \forall
\, \sigma \in {\mathcal{E}}^{int}_h, \;
\alpha_{K_\sigma,L_\sigma}=\frac{1}{2}$.
    Let $K \in {\mathcal{T}}_h$. Using definition (\ref{eq:defg}) and
the gradient formula we get
\begin{align*}
\label{eq:gradu2} \hspace{-.1cm} \int_K \Big( \nabla_h q_h -
\widetilde \nabla_h (\Pi_{P^{nc}_1} q_h) \Big) \, d\x
&=\sum_{\sigma \in {\mathcal{E}}_K \cap {\mathcal{E}}^{int}_h} \frac{|\sigma|}{2} \, (q_{K_\sigma} + q_{L_\sigma}) \, \N_{K,\sigma} \nonumber \\
&+\sum_{\sigma \in {\mathcal{E}}_K \cap
{\mathcal{E}}^{ext}_h}|\sigma| \, q_{K_\sigma} \, \N_{K,\sigma}
 - \sum_{\sigma \in {\mathcal{E}}_K} \int_\sigma (\Pi_{P^{nc}_1} q_h) \, \N_{K,\sigma} \, d\sigma.
\end{align*}
Since $q_h \in P_0$ we deduce from proposition \ref{propp1ncu}
\begin{equation*}
\label{eq:gradu3} \int_\sigma \Pi_{P^{nc}_1} q_h \, d\sigma
=|\sigma| \, (\Pi_{P^{nc}_1} q_h)(\x_\sigma) =\left\{
\begin{array}{lr}
\frac{|\sigma|}{2} \, (q_{K_\sigma} + q_{L_\sigma}) & \hbox{ if } \sigma \in {\mathcal{E}}^{int}_h, \\
|\sigma| \, q_{K_\sigma} & \hbox{ if } \sigma \in
{\mathcal{E}}^{ext}_h.
\end{array}
\right.
\end{equation*}
Plugging this into the equation above, we get $\nabla_h
q_h|_K=\widetilde \nabla_h (\Pi_{P^{nc}_1} q_h)|_K$ . \qed

\begin{lem}
  We assume that the mesh is uniform.  There exists a
  constant $C>0$ such that
  \begin{equation*}\label{eq:infsdeg}
    \forall \, q_h \in P_0 \cap L^2_0 \, , \hspace{1cm} \sup_{\V_h \in \P_0\backslash \{\mathbf{0}\}}
    -\frac{(q_h,\hbox{\rm div}_h \, \V_h)}{\|\V_h\|_h} \ge C \, h \, \|\Pi_{P^{nc}_1} q_h\|_{1,h}.
  \end{equation*}
\end{lem}
\noindent {\sc Proof.} If $q_h=0$ the result is trivial. Let $q_h
\in P_0 \cap L^2_0\backslash \{0\}$. Let $\V_h=\nabla_h q_h \in
\P_0\backslash \{ \mathbf{0}\}$. Using proposition
\ref{prop:propadjh} we have
\begin{equation*}
  -(q_h,\hbox{div}_h \V_h)=(\V_h,\nabla_h q_h)=|\nabla_h q_h|^2=|\nabla_h q_h| \, |\V_h|.
\end{equation*}
Let $\chi_\Omega$ be the characteristic function of $\Omega$. Putting
$\psi=\chi_\Omega$ in (\ref{eq:defpp1nc}) we get $\Pi_{P^{nc}_1} q_h
\in L^2_0$. So according to  (\ref{eq:inpoinp1nc})
and(\ref{eq:gradmu}) we have
\begin{equation*}
  |\nabla_h q_h|=\left|\widetilde \nabla_h(\Pi_{P^{nc}_1} q_h)\right| \ge C \, \|\Pi_{P^{nc}_1} q_h\|_{1,h}.
\end{equation*}
On the other hand, according to proposition \ref{propinvp0}: $
|\V_h| \ge C \, h \, \|\V_h\|_h$. Therefore
\begin{equation*}
 \hspace{3.2cm}  -(q_h,\hbox{div}_h \V_h) \ge C \, h  \, \|\Pi_{P^{nc}_1} q_h\|_{1,h} \, \|\V_h\|_h. \hspace{3.2cm} \qed
\end{equation*}
\begin{prop}
\label{prop:condinfs} We assume that the mesh is uniform. There
exists a constant $C>0$ such that for all $q_h \in P_0 \cap L^2_0$
\begin{equation*}
\label{eq:infsd}
   \sup_{\V_h \in \P_0 \backslash \{\mathbf{0}\}} -\frac{(q_h,\hbox{\rm div}_h \, \V_h)}{\|\V_h\|_h} \ge C
    \, |\Pi_{P^{nc}_1} q_h|.
\end{equation*}
\end{prop}
\noindent {\sc Proof.} If $q_h=0$ the result is clear. Let $q_h \in
P_0 \cap L^2_0 \backslash \{0\}$. According to  (\ref{eq:infsc})
 there exists $\V \in \H^1_0$ such that
\begin{equation}
\label{eq:propv}
  \hbox{div} \, \V=-\Pi_{P^{nc}_1} q_h \hspace{.3cm} \hbox{ and } \hspace{.3cm} \|\V\|_1
   \le C \, |\Pi_{P^{nc}_1} q_h|.
\end{equation}
We set
 $\V_h=\Pi_{\P^c_1} \V$.
We want to  estimate $-\big(q_h,\hbox{div}_h(\Pi_{\P_0} \V_h)\big)$.
Since $\nabla_h q_h \in \P_0$ we deduce from proposition
\ref{prop:propadjh}
\begin{equation*}
\label{eq:infs5}
  -\big(q_h,\hbox{div}_h(\Pi_{\P_0} \V_h)\big)=(\Pi_{\P_0} \V_h,\nabla_h q_h)=(\V_h,\nabla_h q_h).
\end{equation*}
Splitting the last term we get
\begin{equation}
\label{eq:infs4}
  -\big(q_h,\hbox{div}_h(\Pi_{\P_0} \V_h)\big)=(\V,\nabla_h q_h)-(\V-\V_h,\nabla_h q_h).
\end{equation}
One one hand, integrating by parts, we get
\begin{equation*}
(\V,\nabla_h q_h)=-(\Pi_{P^{nc}_1} q_h,\hbox{div} \, \V)+ \sum_{K
\in {\mathcal{T}}_h} \sum_{\sigma \in {\mathcal{E}}_K} \int_{\sigma}
(\Pi_{P^{nc}_1} q_h) \, (\V \cdot \N_{K,\sigma})
 \, d\sigma.
\end{equation*}
According to  (\ref{eq:propv}) we have
  $-(\Pi_{P^{nc}_1} q_h,\hbox{div} \, \V)=|\Pi_{P^{nc}_1} q_h|^2$.
Moreover
\begin{equation*}
\sum_{K \in {\mathcal{T}}_h} \sum_{\sigma \in {\mathcal{E}}_K}
\int_{\sigma} (\Pi_{P^{nc}_1} q_h) \, (\V \cdot \N_{K,\sigma}) \,
d\sigma =\sum_{\sigma \in {\mathcal{E}}^{int}_h} \int_{\sigma}
(\Pi_{P^{nc}_1} q_h) \, (\V \cdot \N_{K_\sigma,\sigma}) \, d\sigma
\end{equation*}
since $\V|_{\partial \Omega}=\mathbf{0}$. Using \cite{brenns} p.269
and (\ref{eq:propv}) we have
\begin{eqnarray*}
 \left|\sum_{\sigma \in {\mathcal{E}}^{int}_h} \int_{\sigma} (\Pi_{P^{nc}_1} q_h)
 \, (\V \cdot \N_{K,\sigma}) \, d\sigma
 \right| &\le& C \, h \, \|\V\|_1 \, \|\Pi_{P^{nc}_1} q_h\|_{1,h} \\
 &\le& C \, h \, |\Pi_{P^{nc}_1} q_h|  \, \|\Pi_{P^{nc}_1} q_h\|_{1,h}.
\end{eqnarray*}
So we get
\begin{equation}
\label{eq:infs6}
  (\V,\nabla_h q_h) \ge (|\Pi_{P^{nc}_1} q_h|-C \, h \, \|\Pi_{P^{nc}_1} q_h\|_{1,h}) \, |\Pi_{P^{nc}_1} q_h|.
\end{equation}
On the other hand, using lemma \ref{lemma:gradmu} and  the
Cauchy-Schwarz inequality
\begin{equation*}
  \left|(\V-\V_h,\nabla_h q_h)\right|=|(\V-\V_h,\widetilde \nabla_h (\Pi_{P^{nc}_1}q_h)
   )| \le |\V-\V_h| \, |\widetilde \nabla_h (\Pi_{P^{nc}_1} q_h)|.
\end{equation*}
Using (\ref{eq:errintp1c}) and (\ref{eq:propv}) we get
\begin{equation*}
  |\V-\V_h|=|\V-\Pi_{\P^c_1} \V| \le C \, h \, \|\V\|_1 \le C \, h \, |\Pi_{P^{nc}_1} q_h|.
\end{equation*}
Thus
\begin{equation*}
\label{eq:infs7}
  \left|(\V-\V_h,\nabla_h q_h)\right| \le C \, h \, |\Pi_{P^{nc}_1} q_h|
   \, |\widetilde \nabla_h(\Pi_{P^{nc}_1} q_h)| \le C \, h \,  |\Pi_{P^{nc}_1} q_h| \, \|\Pi_{P^{nc}_1} q_h\|_{1,h}.
\end{equation*}
Let us plug this  estimate and  (\ref{eq:infs6})  into
(\ref{eq:infs4}).
 We get
\begin{equation*}
\label{eq:infs8}
   -\big(q_h,\hbox{div}_h(\Pi_{\P_0} \V_h)\big) \ge (|\Pi_{P^{nc}_1} q_h|-C \, h \, \|\Pi_{P^{nc}_1}
   q_h\|_{1,h}) \, |\Pi_{P^{nc}_1} q_h|.
\end{equation*}
We now  introduce the norm $\|.\|_h$. We have $\V_h=\Pi_{\P^c_1} \V
\in \P^c_1 \subset \H^1$. Thus, using \cite{eymgal} p. 776, we get
  $\|\Pi_{\P_0} \V_h\|_h \le C \, \|\V_h\|_1$.
Since $\Pi_{\P^c_1}$ is stable for the $\H^1$ norm, we deduce from
(\ref{eq:propv})
\begin{equation*}
  \|\V_h\|_1=\|\Pi_{\P^c_1} \V\|_1 \le \|\V\|_1 \le C \, |\Pi_{P^{nc}_1} q_h|.
\end{equation*}
Therefore
  $\|\Pi_{\P_0} \V_h\|_h \le C \, |\Pi_{P^{nc}_1} q_h|$.
Using this  inequality in (\ref{eq:infs8}) we obtain that there
exists constants  $C_1>0$ and $C_2>0$ such that
\begin{equation*}
   -\big(q_h,\hbox{div}_h(\Pi_{\P_0} \V_h)\big) \ge \left(C_1 \, |\Pi_{P^{nc}_1} q_h|
   -C_2 \, h \, \|\Pi_{P^{nc}_1} q_h\|_{1,h} \right) \,  \|\Pi_{\P_0} \V_h\|_h.
\end{equation*}
We deduce from this
\begin{equation*}
\sup_{\V_h \in \P_0\backslash \{\mathbf{0}\}}
    -\frac{(q_h,\hbox{\rm div}_h \, \V_h)}{\|\V_h\|_h} \ge C_1 \, |\Pi_{P^{nc}_1} q_h|
    -C_2 \, h \, \|\Pi_{P^{nc}_1} q_h\|_{1,h}.
\end{equation*}
Let us combine this  with lemma \ref{eq:infsdeg}. Since
\begin{equation*}
  \forall \, t \ge 0 \, , \hspace{.5cm} \max \big(C \, t \, , \, C_1 \, |\Pi_{P^{nc}_1} q_h| -C_2 \, t \big)
  \ge \frac{C \, C_1}{C+C_2} \, |\Pi_{P^{nc}_1} q_h|
\end{equation*}
we get the result. \qed

\subsection{Properties of the discrete laplacian}
\label{subsec:proplapd} We first prove the coercivity of the
discrete laplacian.

\begin{prop}
\label{prop:coerlap}

\noindent  For all $\U_h \in \P_0$ and $\V_h \in \P_0$ we have
  \begin{equation*}
  -(\Deltat_h \U_h,\U_h)=\|\U_h\|^2_h \hspace{1cm}
    -(\Deltat_h \U_h,\V_h) \le \|\U_h\|_h \, \|\V_h\|_h.
  \end{equation*}
\end{prop}
\noindent {\sc Proof.} Using definition (\ref{eq:defl}) and writing
 the sum over the triangles as a sum over the edges, we have
\begin{eqnarray*}
\label{eq:deltnh}
  -(\Deltat_h \U_h, \V_h) &=&
 -\sum_{K \in {\mathcal{T}}_h} \V_K \cdot  \Big( \sum_{\sigma \in {\mathcal{E}}_K \cap
{\mathcal{E}}^{int}_h}
  \tau_{\sigma} \, (\U_{L_\sigma}-\U_{K})      -
\sum_{\sigma \in {\mathcal{E}}_K \cap {\mathcal{E}}^{ext}_h} \tau_\sigma \, \U_K  \Big) \nonumber \\
 &=& \sum_{\sigma \in {\mathcal{E}}^{int}_h} \tau_{\sigma} \,
(\V_{L_\sigma}-\V_K) \cdot (\U_{L_\sigma}-\U_K) +\sum_{K \in
{\mathcal{T}}_h} \sum_{\sigma \in {\mathcal{E}}_K \cap
{\mathcal{E}}^{ext}_h} \tau_\sigma \, \U_K \cdot \V_K.
\end{eqnarray*}
We get the first half of the result by taking $\V_h=\U_h$. On the
other hand, using the Cauchy-Schwarz inequality and  the algebraic
identity $a \, b+c \, d \le \sqrt{a^2+c^2} \sqrt{b^2+d^2}$, we get
the second half.\qed


\noindent  If $\V \in \H^2$, we have  $|\Delta \V| \le \|\V\|_2$.
The operator
 $\mathbf{\Delta}_h$ shares a similar property.
\begin{prop}
\label{prop:stablap}
  There exists a constant $C>0$ such that for all  $\V \in \H^2$
  \begin{equation*}
    \left|\Deltat_h(\widetilde \Pi_{\P_0} \V)\right| \le C \, \|\V\|_2.
  \end{equation*}
\end{prop}
\noindent {\sc Proof.} Let $\V_h=\widetilde \Pi_{\P_0} \V$. Let $K
\in {\mathcal{T}}_h$. According to definition (\ref{eq:deflap})
\begin{equation}
\label{eq:stablap1}
  \Deltat_h \V_h|_K=\frac{1}{|K|} \sum_{\sigma \in {\mathcal{E}}_K \cap {\mathcal{E}}^{int}_h} \tau_\sigma \,
  (\V(\x_{L_\sigma})-\V(\x_K)) - \frac{1}{|K|} \sum_{\sigma \in {\mathcal{E}}_K \cap {\mathcal{E}}^{ext}_h}
   \tau_\sigma \, \V(\x_K).
\end{equation}
Let us first assume that $\V=(v_1,v_2) \in
({\mathcal{C}}^\infty_0)^2$. Let $i\in\{1,2\}$. If $\sigma \in
{\mathcal{E}}_K \cap {\mathcal{E}}^{int}_h $ and $\x \in \sigma$ we
have the Taylor expansions
\begin{equation*}
\label{eq:stablap2}
  v_i(\x_{L_\sigma})=v_i(\x)+\nabla v_i(\x) \cdot (\x_{L_\sigma}-\x)
  +\int^1_0 \H(v_i)(t \, \x_{L_\sigma}+(1-t) \, \x) (\x_{L_\sigma}-\x) \cdot (\x_{L_\sigma}-\x)\, t \, dt \, ,
\end{equation*}
\begin{equation*}
\label{eq:stablap3}
  v_i(\x_{K})=v_i(\x)+\nabla v_i(\x) \cdot (\x_{K}-\x)
  +\int^1_0 \H(v_i)(t \, \x_{K}+(1-t) \, \x) (\x_{K}-\x) \cdot (\x_{K}-\x)\, t \, dt \, ,
\end{equation*}
\begin{equation*}
\label{eq:stablap4}
  \nabla v_i(\x)=\nabla v_i(\x_K)-\int^1_0 \Nabla \nabla v_i(t \, \x_K+(1-t) \, \x)(\x_K-\x) \, dt.
\end{equation*}
The notation $\H(v_i)$ refers to the hessian matrix of $v_i$.
Plugging the last expansion into the two others and integrating over
$\sigma$, we get
\begin{equation*}
\label{eq:stablap5}
  \int_\sigma \big(v_i(\x_{L_\sigma})-v_i(\x)\big) \, d\x=\nabla v_i(\x_K) \cdot (\x_{L_\sigma}-\x_\sigma)
  -A^{\sigma,i}_{L_\sigma}+B^{\sigma,i}_{L_\sigma} \, ,
\end{equation*}
\begin{equation*}
\label{eq:stablap6}
  \int_\sigma \big(v_i(\x_{K})-v_i(\x)\big) \, d\x=\nabla v_i(\x_K) \cdot (\x_{K}-\x_\sigma)
  -A^{\sigma,i}_{K}+B^{\sigma,i}_{K}.
\end{equation*}
The terms $A^{\sigma,i}_T$ and $B^{\sigma,i}_T$ are the same as in
(\ref{eq:defah}) and (\ref{eq:defbh}), with $v_i$ instead of $q$.
We substract these  equations. Since
$\x_{L_\sigma}-\x_K=d_\sigma \, \N_{K,\sigma}$ we infer from
(\ref{eq:deftaus})
\begin{equation*}
\label{eq:stablap65}
 \hspace{.2cm} \tau_\sigma \, \big(v_i(\x_{L_\sigma})-v_i(\x_K)\big)= \nabla v_i(\x_K)\cdot
  \N_{K,\sigma}
  +\frac{1}{d_\sigma} \,
  \left(-A^{\sigma,i}_{L_\sigma}+B^{\sigma,i}_{L_\sigma}
+A^{\sigma,i}_{K}-B^{\sigma,i}_{K}\right).
\end{equation*}
Let us consider now the case $\sigma \in {\mathcal{E}}_{K} \cap
{\mathcal{E}}^{ext}_h$. If $\x \in \sigma$ we have the Taylor
expansions
\begin{equation*}
\label{eq:stablap7}
 v_i(\x_{K})=v_i(\x)+\nabla v_i(\x) \cdot (\x_{K}-\x)
  +\int^1_0 \H(v_i)(t \, \x_{K}+(1-t) \, \x) (\x_{K}-\x) \cdot (\x_{K}-\x)\, t \, dt \, ,
\end{equation*}
\begin{equation*}
\label{eq:stablap8}
    \nabla v_i(\x)=\nabla v_i(\x_K)-\int^1_0 \Nabla \nabla v_i(t \, \x_K+(1-t) \, \x)(\x_K-\x) \, dt.
\end{equation*}
Since $v_i \in {\mathcal{C}}^\infty_0$ we have $v_i(\x)=0$. We
plug the last expansion  into the other and integrate over $\sigma$.
Since $\x_K-\x_\sigma=-d_\sigma \, \N_{K,\sigma}$ we deduce from
(\ref{eq:deftaus})
\begin{equation*}
\label{eq:stablap75}
  -\tau_\sigma \, v_i(\x_K)= \nabla v_i(\x_K) \cdot \N_{K,\sigma}
  +\frac{1}{d_\sigma} \,  \big(A^{\sigma,i}_K-B^{\sigma,i}_K\big).
\end{equation*}
Thus we get
\begin{align*}
  \frac{1}{|K|} \sum_{\sigma \in {\mathcal{E}}_K \cap
{\mathcal{E}}^{int}_h} \tau_\sigma \,
\big(v_i(\x_{L_\sigma})&-v_i(\x_K)\big)
- \frac{1}{|K|} \sum_{\sigma \in {\mathcal{E}}_K \cap {\mathcal{E}}^{ext}_h} \tau_\sigma \, v_i(\x_K) \\
&=\frac{1}{|K|} \, \nabla v_i(\x_K) \cdot \sum_{\sigma \in
{\mathcal{E}}_K} |\sigma| \, \N_{K,\sigma}+ \frac{1}{|K|} \,
\sum_{\sigma \in {\mathcal{E}}_K} R^i_\sigma
\end{align*}
where we have set for all edge $\sigma \in {\mathcal{E}}_K \cap
{\mathcal{E}}^{int}_h$
\begin{equation*}
  R^i_\sigma=\frac{1}{d_\sigma} \, \left(-A^{\sigma,i}_{L_\sigma}+B^{\sigma,i}_{L_\sigma}
+A^{\sigma,i}_{K}-B^{\sigma,i}_{K}\right)
\end{equation*}
and for all edge  $\sigma \in {\mathcal{E}}_K \cap
{\mathcal{E}}^{ext}_h$:
  $R^i_\sigma=\frac{1}{d_\sigma} \,  \big(A^{\sigma,i}_K-B^{\sigma,i}_K\big)$.
Since $\sum_{\sigma \in {\mathcal{E}}_K} |\sigma| \,
\N_{K,\sigma}=\mathbf{0}$, setting
$\mathbf{R}_\sigma=(R^1_\sigma,R^2_\sigma)$, we get
\begin{equation*}
 \frac{1}{|K|} \sum_{\sigma \in {\mathcal{E}}_K \cap {\mathcal{E}}^{int}_h} \tau_\sigma \, \big(\V(\x_{L_\sigma})-\V(\x_K)\big)
- \frac{1}{|K|} \sum_{\sigma \in {\mathcal{E}}_K \cap
{\mathcal{E}}^{ext}_h} \tau_\sigma \, \V(\x_K) = \frac{1}{|K|}
\sum_{\sigma \in {\mathcal{E}}_K} \R_\sigma.
\end{equation*}
Since the space  $({\mathcal{C}}^\infty_0)^2$ is dense in $\H^2$, one checks that
this equation still holds for  $\V \in \H^2$.
Using (\ref{eq:stablap1}) we infer from it
\begin{equation*}
  \left| \Deltat_h \V_h \right|^2=\sum_{K \in {\mathcal{T}}_h}  |K|\left| \Deltat_h \V_h|_K \right|^2
  \le \sum_{K \in {\mathcal{T}}_h} \frac{3}{|K|}\sum_{\sigma \in {\mathcal{E}}_K} |\R_\sigma|^2.
\end{equation*}
Using estimates (\ref{eq:defah}) and (\ref{eq:defbh}) we obtain
\begin{equation*}
\hspace{1.4cm}  \left| \Deltat_h \V_h \right|^2 \le
  C  \, \sum_{i=1}^2 \sum_{\sigma \in {\mathcal{E}}_h} \int_{D_\sigma}
   \left( \left|\Nabla \nabla v_i\right|^2+\left|\H(v_i)\right|^2 \right)\, d\x \le C \, \|\V\|^2_2.
\hspace{1.4cm}    \qed
\end{equation*}

\section{Stability of the scheme}
\label{sec:stab}

We now use the results of section \ref{sec:propopd} to prove the
stability of the scheme.  We first show an estimate for the computed
velocity (theorem \ref{theo:estv}).
 We then state a similar result for the increments in time (lemma
 \ref{lem:estiv}).
Using the inf-sup condition (proposition \ref{prop:condinfs}),  we infer from it  some estimates on the pressure (theorem \ref{theo:estp}).


\begin{lem}
\label{lem:orthp} For all $m\in\{0,\dots,N\}$ and
$n\in\{0,\dots,N\}$ we have
\begin{eqnarray*}
     (\U^m_h,\nabla_h p^n_h)=0  \, , \hspace{1cm}
     |\U^m_h|^2-|\Ut^m_h|^2+|\U^m_h-\Ut^m_h|^2=0. \label{eq:orth2}
\end{eqnarray*}
\end{lem}
\noindent {\sc Proof.} First, using propositions \ref{prop:umrt0}
and \ref{prop:propadjh}, we get
\begin{equation*}
 (\U^m_h,\nabla_h p^n_h)=-(p^n_h,\hbox{div}_h \U^m_h)=0.
\end{equation*}
Thus we deduce from   (\ref{eq:projib})
\begin{equation*}
 2 \, (\U^m_h,\U^m_h-\Ut^m_h)
=-\frac{4 \, k}{3} \big(\U^m_h,\nabla_h(p^m_h-p^{m-1}_h)\big)=0.
\end{equation*}
Using the  algebraic identity $2 \, a \, (a-b)=a^2-b^2+(a-b)^2$ we
get
\begin{equation*}
\hspace{2.1cm} 2 \,
(\U^m_h,\U^m_h-\Ut^m_h)=|\U^m_h|^2-|\Ut^m_h|^2+|\U^m_h-\Ut^m_h|^2=0.
\hspace{2.1cm} \qed
\end{equation*}
We introduce the following hypothesis on the initial data.
\begin{equation*}
  {\bf (H1)}  \hspace{.2cm} \hbox{ There exists $C>0$ such that }
    \hspace{.2cm} |\U^0_h|+|\U^1_h|+k|\nabla_h p^1_h| \le C.
\end{equation*}
Hypothesis {\bf (H1)} is fulfilled if we set
$\U^0_h=\Pi_{\mathbf{RT_0}} \U_0$ and we use a semi-implicit Euler
 scheme to compute  $\U^1_h$. We have the following result.
\begin{theo}
\label{theo:estv} We assume  that the initial values  of the scheme
fulfill {\bf (H1)}.  For all $m \in \{2,\dots,N\}$ we have
\begin{equation*}
\label{eq:estl2h1}
  |\U^m_h|^2+k\sum_{n=2}^m \|\Ut^n_h\|^2_h  \le C.
\end{equation*}
\end{theo}
\noindent {\sc Proof.} Let $m \in\{2,\dots,N\}$
 and $n \in \{1,\dots,m-1\}$. Taking the scalar product of
 (\ref{eq:mombdf}) with
$4 \, k \, \Ut^{n+1}_h$ we get
\begin{align}
\label{eq:eqstabmultbdf}
    \left( \frac{3 \, \Ut^{n+1}_h-4 \, \U^n_h+\U^{n-1}_h}{2k}, 4 \, k \, \Ut^{n+1}_h \right)
&-\frac{4 \, k}{\hbox{Re}} \, (\Deltat_h \Ut^{n+1}_h, \Ut^{n+1}_h)
  \nonumber \\
 +4 \, k \, \b_h(2 \, \U^n_h -\U^{n-1}_h,\Ut^{n+1}_h,\Ut^{n+1}_h) &+4 \, k \, (\nabla_h p^n_h,\Ut^{n+1}_h)
=4 \, k \, (\F^{n+1}_h,\Ut^{n+1}_h).
\end{align}
First of all, using lemma \ref{lem:orthp}, we get as in \cite{guer2}
\begin{eqnarray*}
\label{eq:estvbdf1}
&& 4\, k \left(\Ut^{n+1}_h,\frac{3 \, \Ut^{n+1}_h-4 \, \U^n_h+\U^{n-1}_h}{2 \, k}\right) \nonumber \\
 &&= |\U^{n+1}_h|^2-|\U^n_h|^2 +6 \, |\Ut^{n+1}_h-\U^{n+1}_h|^2
    +|2 \, \U^{n+1}_h - \U^n_h|^2
  -|2\, \U^n_h - \U^{n-1}_h|^2 \nonumber \\
  &&+|\U^{n+1}_h -2 \, \U^n_h +\U^{n-1}_h|^2.
\end{eqnarray*}
According to  proposition \ref{prop:coerlap} we have
$ -\frac{4 \, k}{\hbox{Re}} \, (\Deltat_h
\Ut^{n+1}_h, \Ut^{n+1}_h)=\frac{4 \, k}{\hbox{Re}} \,
\|\Ut^{n+1}_h\|^2_h$.
Also, using  lemma \ref{lem:orthp} and
(\ref{eq:projib}), we have
\begin{eqnarray*}
\label{eq:estvbdf4}
4 \, k \, (\nabla_h p^n_h,\Ut^{n+1}_h)&=&4 \, k \, (\nabla_h p^n_h,\Ut^{n+1}_h-\U^{n+1}_h) \nonumber \\
&=&\frac{4 \, k^2}{3} \, ( |\nabla  p^{n+1}_h|^2 -
 |\nabla p^{n}_h|^2 - |\nabla  p^{n+1}_h -\nabla p^{n}_h |^2 ).
\end{eqnarray*}
Multiplying  (\ref{eq:projib}) by $4k\nabla_h (p^{n+1}_h-p^n_h)$ and
using the Young inequality we get
\begin{equation*}
\label{eq:estvbdf5} \frac{4 \, k^2}{3} \, |\nabla
(p^{n+1}_h-p^n_h)|^2 \le 3 \, |\U^{n+1}_h - \Ut^{n+1}_h|^2.
\end{equation*}
According to proposition \ref{prop:posbh} we have
  $4 \, k \, \b_h(2 \, \U^n_h - \U^{n-1}_h,\Ut^{n+1}_h,\Ut^{n+1}_h) \ge 0$.
At last using the  Cauchy-Schwarz inequality, (\ref{eq:inpoinp0})
and (\ref{eq:propfh}) we have
\begin{equation*}
 4 \, k \, (\F^{n+1}_h, \Ut^{n+1}_h) \le 4 \, k \, |\F^{n+1}_h| \, |\Ut^{n+1}_h|
  \le C \, k \, \|\F\|_{{\mathcal{C}}(0,T;\L^2)}\, \|\Ut^{n+1}_h\|_h.
 \end{equation*}
Using the  Young inequality we get
\begin{equation*}
\label{eq:estvbdf3}
4 \, k \, (\F^{n+1}_h, \Ut^{n+1}_h) \le  3 \, k
\, \|\Ut^{n+1}_h\|^2_h + C \, k \,
\|\F\|^2_{{\mathcal{C}}(0,T;\L^2)}.
\end{equation*}
Let us plug these estimates   into (\ref{eq:eqstabmultbdf}). We get
\begin{eqnarray*}
\label{eq:estvbdf6}
 \hspace{-1cm} & & |\U^{n+1}_{h}|^2 - |\U^n_{h}|^2
+ |2 \, \U^{n+1}_h-\U^n_h|^2-|2 \, \U^n_h-\U^{n-1}_h|^2
 +|\U^{n+1}_h-2 \, \U^n_h+\U^{n-1}_h|^2 \nonumber \\
\hspace{-1cm} && +3 \, |\Ut^{n+1}_h - \U^{n+1}_h|^2
 +k \, \|\Ut^{n+1}_{h}\|^2_h +\frac{4 \, k^2}{3} \, (|\nabla_h p^{n+1}_h|^2 - |\nabla_h p^n_h|^2)  \le C \, k.
\end{eqnarray*}
Summing from $n=1$ to $m-1$ we have
\begin{eqnarray*}
\label{eq:estvbdf7} \hspace{-.6cm} && |\U^m_h|^2+|2 \,
\U^m_h-\U^{m-1}_h|^2 +3\sum_{n=1}^{m-1} |\Ut^{n+1}_h - \U^{n+1}_h|^2
+k\sum_{n=1}^{m-1} \|\Ut^{n+1}_h\|^2_h + \frac{4 \, k^2}{3} \, |\nabla_h p^m_h|^2  \nonumber \\
\hspace{-.6cm}&& \le C+4 \, |\U^1_h|^2+|2 \, \U^1_h - \U^0_h|^2+k^2
\, |\nabla_h p^1_h|^2.
\end{eqnarray*}
Using hypothesis {\bf (H1)} we get the result. \qed

We now want to estimate the computed pressure.
From now on, we make the following hypothesis on the data
\begin{equation*}
    \F \in {\mathcal{C}}(0,T;\L^2) \, , \hspace{.6cm} \F_t \in L^2(0,T;\L^2) \, ,
  \hspace{.6cm} \U_0 \in \H^2 \cap \H^1_0 \, , \hspace{.6cm} \hbox{div} \, \U_0=0.
\end{equation*}
 For all sequence  $(q^m)_{m \in
\mathbb{N}}$ we define the sequence  $(\delta q^m)_{m \in
\mathbb{N}^*}$ by setting $\delta
q^m=q^m-q^{m-1}$ for $m \ge 1$.
 We set $\Del=(\delta)^2$. If the data
$\U_0$ and $\F$ fulfill
 a compatibility condition  \cite{heywood}
there exists a  solution  $(\U,p)$ to the  equations
(\ref{eq:mom})--(\ref{eq:incomp}) such that
\begin{equation*}
\label{eq:regup}
 \U \in {\mathcal{C}}(0,T;\H^2) \, , \hspace{.5cm} \U_t \in {\mathcal{C}}(0,T;\L^2) \, ,
  \hspace{.5cm} \nabla p\in {\mathcal{C}}(0,T;\L^2).
\end{equation*}
We introduce the following hypothesis on the initial values of the
scheme: there exists a constant $C>0$ such that
\begin{equation*}
  {\bf (H2)}    \hspace{.7cm}  |\U^0_h-\U_0| + \frac{1}{h} \, \|\U^1_h-\U(t_1)\|_{\infty}
      +|p^1_h-p(t_1)| \le C \, h \, , \hspace{.7cm} |\U^1_h-\U^0_h|\le C \, k.
\end{equation*}
One checks easily that this hypothesis implies {\bf (H1)}. We have
the following result.
\begin{lem}
\label{lem:estiv} We assume that the initial values of the scheme
fulfill {\bf (H2)}. Then there exists a constant $C>0$ such that for
all $m \in \{1,\dots,N\}$
\begin{equation}
\label{eq:estincv}
  \frac{1}{k} \, |\Del \U^m_h| \le C .
\end{equation}
\end{lem}
\noindent {\sc Proof.} We prove the result by induction. The result
holds for  $m=1$ 
thanks to hypothesis {\bf (H2)}.
 Let us consider the case  $m=2$.
We  set $\Ut^1_h=\U^1_h$. Let $\U^{-1}_h \in \P_0$ given by
\begin{equation}
\label{eq:estiv89}
  \U^{-1}_h=4 \, \U^0_h-3 \, \U^1_h+\frac{2 \, k}{\hbox{Re}} \,  \Deltat_h \U^1_h
-2 \, k\, \bt_h(\U^0_h,\Ut^1_h)  -2 \, k \, \nabla_h p^1_h -2 \, k
\, \F^1_h.
\end{equation}
We substract this equation from equation (\ref{eq:projib}) written
for $n=1$. Since
\begin{equation*}
\bt_h(2 \, \U^1_h-\U^0_h,\Ut^2_h)-\bt_h(\U^0_h,\Ut^1_h) = \bt(2 \,
\U^1_h-2 \, \U^0_h,\Ut^2_h)+\bt_h(\U^0_h,\Del \Ut^2_h) \, ,
\end{equation*}
upon setting  $\Del \U^0_h=\U^0_h-\U^{-1}_h$, we get
\begin{equation*}
\label{eq:estiv9}
  \frac{3 \, \Del \Ut^2_h - 4 \, \Del \U^1_h+\Del \U^0_h}{2 \, k} - \frac{1}{\hbox{Re}} \, \Deltat_h (\Del \Ut^2_h)
  +\bt_h(2 \, \U^1_h-2 \, \U^0_h,\Ut^2_h)+\bt_h(\U^0_h,\Del \Ut^2_h)=\Del \F^2_h.
\end{equation*}
Taking the scalar product with  $4 \, k \, \Del \Ut^2_h$ we get
\begin{eqnarray}
\label{eq:estiv10} \hspace{-1cm} && 2 \left(  3 \, \Del \Ut^2_h - 4
\, \Del \U^1_h+\Del \U^0_h, \Del \Ut^2_h \right) -
\frac{1}{\hbox{Re}} \, \big( \Deltat_h (\Del \Ut^2_h),\Del \Ut^2_h
\big)
 \nonumber \\
&&+ 4 \, k \, \b_h(\U^0_h,\Del \Ut^2_h,\Del \Ut^2_h)+4 \, k \,
\b_h(2 \, \U^1_h-2 \, \U^0_h,\Ut^2_h,\Del \Ut^2_h) =4 \, k\, (\Del
\F^2_h,\Del \Ut^2_h).
\end{eqnarray}
According to  proposition \ref{prop:stabbth} we have
\begin{equation*}
  4 \, k \, |\b_h(2 \, \U^1_h-2 \, \U^0_h,\Ut^2_h,\Del \Ut^2_h)|
  \le C \, k \, |2 \, \U^1_h - 2 \, \U^0_h| \, \|\Ut^2_h\|_h \, \|\Del \Ut^2_h\|_h \, ;
\end{equation*}
so that, using  hypothesis {\bf (H2)}
\begin{equation*}
  4 \, k \, \left|\b_h(2 \, \U^1_h-2 \, \U^0_h,\Ut^2_h,\Del \Ut^2_h)\right|
  \le C \, k^2 \, \|\Ut^2_h\|_h \, \|\Del \Ut^2_h\|_h.
\end{equation*}
From the  Young inequality and  theorem \ref{theo:estv} we deduce
\begin{equation*}
 4 \, k \, \left|\b_h(2 \, \U^1_h-\U^0_h,\Ut^2_h,\Ut^2_h-\Ut^1_h)\right|
 \le \frac{k}{\hbox{Re}} \, \|\Del \Ut^2_h\|^2_h+C \, k^3 \, \|\Ut^2_h\|^2_h
 \le  \frac{k}{\hbox{Re}} \, \|\Del \Ut^2_h\|^2_h+C \, k^2.
\end{equation*}
On the other hand
\begin{equation*}
  \Del \F^2_h= \F^2_h - \F^1_h=\Pi_{\P_0} \F(t_2)-\Pi_{\P_0}\F(t_1)=\Pi_{\P_0} \left( \int^{t_2}_{t_1} \F_t(s) \, ds \right).
\end{equation*}
Since $\Pi_{\P_0}$ is stable for the $\L^2$ norm, using the
Cauchy-Schwarz inequality, we get
\begin{equation*}
  |\Del \F^2_h| \le \int^{t_2}_{t_1} |\F_t(s)| \, ds \le
\sqrt{k} \, \left( \int^{t_2}_{t_1} |\F_t(s)|^2 \, ds \right)^{1/2}
  \le  \sqrt{k} \, \|\F_t\|_{L^2(0,T;\L^2)}.
\end{equation*}
Thus
\begin{equation*}
  4 \, k \, |(\Del \F^2_h,\Del \Ut^2_h)| \le 4 \, k \, |\Del \F^2_h| \, |\Del \Ut^2_h|
  \le C \, k^{3/2} \, |\Del \Ut^2_h|.
\end{equation*}
So that, using (\ref{eq:inpoinp0}) and the  Young inequality
\begin{equation*}
 4 \, k \, |(\Del \F^2_h,\Del \Ut^2_h)|  \le C \, k^{3/2} \, \|\Del \Ut^2_h\|_h
 \le \frac{k}{\hbox{Re}} \, \|\Del \Ut^2_h\|^2_h +C \, k^2.
\end{equation*}
The other terms in  (\ref{eq:estiv10}) are dealt with as in the
prooof of theroem  \ref{theo:estv}. We get
\begin{equation}
\label{eq:estiv104}
  |\Del \U^2_h|^2 \le |\Del \U^1_h|^2+|2 \, \Del \U^1_h-\Del \U^0_h|^2.
\end{equation}
We know ((\ref{eq:estincv}) for $m=1$) that  $|\Del \U^1_h|^2 \le C
\, k^2$.  It remains to estimate the term  $|2 \, \Del \U^1_h-\Del
\U^0_h|^2$. According to  (\ref{eq:estiv89})
\begin{equation*}
  2 \, \Del \U^1_h-\Del \U^0_h=-\Del \U^1_h+\frac{2 \, k}{\hbox{Re}} \, \Deltat_h \U^1_h-2 \, k \, \bt_h(\U^0_h,\U^1_h)
  -2 \, k \, \nabla_h p^1_h -2 \, k \, \F^1_h \, ;
\end{equation*}
by taking the scalar product with  $2 \, \Del \U^1_h-\U^0_h$ and
using the Cauchy-Schwarz inequality we get
\begin{eqnarray}
\label{eq:estiv105}
 |2 \, \Del \U^1_h-\Del \U^0_h|^2 &\le& 2 \, k \, \Big( \frac{|\Del \U^1_h|}{2 \, k} +\frac{1}{\hbox{Re}} \, |\Deltat_h \U^1_h|
  +|\nabla_h p^1_h| +|\F^1_h|\Big) |2 \, \Del \U^1_h -\Del \U^0_h| \nonumber \\
 &+&2 \, k \, \left|\b(\U^0_h,\Ut^1_h,2 \, \Del \U^1_h-\Del \U^0_h)\right|.
\end{eqnarray}
Let us bound the terms between braces. First, we have
\begin{equation*}
\label{eq:estiv11}
  \Deltat_h \U^1_h=\Deltat_h \big(\U^1_h-\widetilde \Pi_{\P_0} \U(t_1)\big)+
\Deltat_h  \big(\Pi_{\P_0} \U(t_1)\big).
\end{equation*}
On one hand, according to  proposition \ref{prop:coerlap}
\begin{eqnarray*}
\left|\Deltat_h \big(\U^1_h-\widetilde \Pi_{\P_0}
\U(t_1)\big)\right|^2
&=&\left(\Deltat_h \big(\U^1_h-\widetilde \Pi_{\P_0} \U(t_1)\big),\Deltat_h \big(\U^1_h-\widetilde \Pi_{\P_0} \U(t_1)\big)\right) \\
& \le & \| \Deltat_h \big(\U^1_h-\widetilde \Pi_{\P_0} \U(t_1)\big)
\|_h \, \|\U^1_h-\widetilde \Pi_{\P_0} \U(t_1)\|_h.
\end{eqnarray*}
Applying  proposition \ref{propinvp0} we get
\begin{equation*}
\left|\Deltat_h \big(\U^1_h-\widetilde \Pi_{\P_0}
\U(t_1)\big)\right|^2 \le \frac{C}{h^2} \, |\Deltat_h
\big(\U^1_h-\widetilde \Pi_{\P_0} \U(t_1)\big)| \,
|\U^1_h-\widetilde \Pi_{\P_0} \U(t_1)|.
\end{equation*}
Using the embedding  $\L^\infty \subset \L^2$  we have
\begin{equation*}
  |\U^1_h-\widetilde \Pi_{\P_0} \U(t_1)|=|\widetilde \Pi_{\P_0}( \U^1_h-\U(t_1))|
  \le \|\widetilde \Pi_{\P_0}( \U^1_h-\U(t_1))\|_\infty \, ;
\end{equation*}
since $\widetilde \Pi_{\P_0}$ is stable for the $\L^\infty$ norm, we
get  using hypothesis {\bf (H2)}
\begin{equation*}
\label{eq:estiv119}
   |\U^1_h-\widetilde \Pi_{\P_0} \U(t_1)| \le   \| \U^1_h-\U(t_1)\|_\infty \le C \, h^2.
\end{equation*}
Therefore
  $\left|\Deltat_h \big(\U^1_h-\widetilde \Pi_{\P_0} \U(t_1)\big)\right| \le C$.
And according to  proposition \ref{prop:stablap}
\begin{equation*}
  \left|\Deltat_h  \big(\Pi_{\P_0} \U(t_1)\big)\right| \le C \, \|\U(t_1)\| \le C\, \|\U\|_{{\mathcal{C}}(0,T;\H^2)}.
\end{equation*}
Hence 
$|\Deltat_h \U^1_h| \le C$.
Let us now bound the pressure term in  (\ref{eq:estiv105}). We have
\begin{equation*}
  \nabla_h p^1_h= \nabla_h\big(p^1_h - \widetilde \Pi_{P_0} p(t_1)\big) +
  \left( \nabla_h \big(\widetilde \Pi_{P_0} p(t_1)\big) - \Pi_{\P_0} \nabla p(t_1)\right)
  +\Pi_{\P_0} \nabla p(t_1).
\end{equation*}
According to proposition \ref{prop:ininvgrad} we have
$\left|\nabla_h\big(p^1_h - \widetilde \Pi_{P_0} p(t_1)\big)\right|
\le \frac{C}{h} \, |p^1_h - \widetilde \Pi_{P_0} p(t_1)|$.
Using (\ref{eq:estdiffprojp0}) we get
\begin{equation*}
\left|\nabla_h\big(p^1_h - \widetilde \Pi_{P_0} p(t_1)\big)\right|
 \le C  \, \|p(t_1)\|_2 \le
C \, \|p\|_{{\mathcal{C}}(0,T;\H^2)}.
\end{equation*}
Since $\P_0$ is stable for the $\L^2$ norm we have
$ \left|\Pi_{\P_0}\nabla p(t_1)\right|
 \le  \left|\nabla p(t_1)\right| \le  \|p\|_{{\mathcal{C}}(0,T;\H^1)}$.
Using  proposition \ref{prop:consg} to treat last term we get
  $|\nabla_h p^1_h| \le C$.
And according to  (\ref{eq:propfh}) and (\ref{eq:estincv}) for $m=1$
we have
  $\frac{|\Del \U^1_h|}{2 \, k}+|\F^1_h| \le C$.
We are left with the term  $\left|\b_h(\U^0_h,\Ut^1_h,2 \, \Del
\U^1_h-\Del \U^0_h)\right|$ in (\ref{eq:estiv105}). We use the
following splitting
\begin{eqnarray*}
  \bt_h(\U^0_h,\U^1_h)&=&\bt_h(\U^0_h-\Pi_{\mathbf{RT_0}} \U_0,\U^1_h)
  +\bt_h\big(\Pi_{\mathbf{RT_0}} \U_0,\U^1_h-\widetilde \Pi_{\P_0}
  \U(t_1)\big) \\
  &+&\bt_h\big(\Pi_{\mathbf{RT_0}}\U_0,\widetilde \Pi_{\P_0} \U(t_1)\big).
\end{eqnarray*}
Let us take the scalar product  with  $2 \, \Del \U^1_h-\Del
\U^0_h$. We get
\begin{equation*}
\label{eq:estiv15} \b_h(\U^0_h,\U^1_h,2 \, \Del \U^1_h -\Del
\U^0_h)=B_1+B_2+B_3
\end{equation*}
with
\begin{equation*}
  B_1=\b_h(\U^0_h-\Pi_{\mathbf{RT_0}} \U_0,\U^1_h,2 \, \Del \U^1_h-\Del
  \U^0_h)\, ,
\end{equation*}
\begin{equation*}
B_2= \b_h\big(\Pi_{\mathbf{RT_0}} \U_0,\U^1_h-\widetilde \Pi_{\P_0}
\U(t_1),2 \, \Del \U^1_h-\Del \U^0_h\big) \, ,
\end{equation*}
and
\begin{equation*}
  B_3=\Big(\bt_h\big(\Pi_{\mathbf{RT_0}}\U_0,\widetilde \Pi_{\P_0} \U(t_1)\big),2 \, \Del \U^1_h-
  \Del\U^0_h\Big).
\end{equation*}
Applying propositions \ref{propinvp0} and \ref{prop:stabbth} we have
\begin{equation*}
|B_1|
  \le \frac{C}{h} \, |\U^0_h-\Pi_{\mathbf{RT_0}} \U_0| \, \|\U^1_h\|_h \, |2 \, \Del \U^1_h-\Del \U^0_h|.
\end{equation*}
According to  (\ref{eq:estdiffprojp0}) and (\ref{eq:esterrirt0})  we have
have
\begin{equation*}
|\U^0_h-\Pi_{\mathbf{RT_0}} \U_0|=
|\Pi_{\P_0}\U^0-\Pi_{\mathbf{RT_0}} \U_0|\le |\Pi_{\P_0}
\U_0-\U_0|+|\U_0-\Pi_{\mathbf{RT_0}} \U_0| \le C \, h \, \|\U_0\|_1.
\end{equation*}
According to  proposition \ref{prop:coerlap} and (\ref{eq:inpoinp0})
\begin{equation*}
  \|\U^1_h\|^2_h=-(\Deltat_h \U^1_h,\U^1_h) \le |\Deltat_h \U^1_h| \, |\U^1_h| \le C \, |\Deltat_h \U^1_h| \, \|\U^1_h\|_h \, ;
\end{equation*}
since $|\Deltat_h \U^1_h|$ is bounded we get
 $\|\U^1_h\|_h \le C$.
Hence
   $|B_1| \le C \, |2 \, \Del \U^1_h -\Del \U^0_h|$.
In a similar way, using propositions \ref{propinvp0} and
\ref{prop:stabbth}, we get
\begin{equation*}
|B_2| \le \frac{C}{h^2} \, |\Pi_{\mathbf{RT_0}} \U_0| \,
|\U^1_h-\widetilde \Pi_{\P_0} \U(t_1)| \,|2 \, \Del \U^1_h-\Del
\U^0_h|.
\end{equation*}
We have
  $|\Pi_{\mathbf{RT_0}} \U_0| \le |\Pi_{\mathbf{RT_0}} \U_0-\U_0|+|\U_0| \le C \, h \, \|\U_0\|_1+|\U_0| \le C \,
  \|\U_0\|_1$.
Using moreover (\ref{eq:estiv119}) we get
$|B_2| \le C \, |2 \, \Del \U^1_h-\Del \U^0_h|$.
Lastly  using the following splitting
\begin{eqnarray*}
\bt_h\big(\Pi_{\mathbf{RT_0}}\U_0,\widetilde \Pi_{\P_0} \U(t_1)\big)
&=&\Big( \bt_h\big(\Pi_{\mathbf{RT_0}}\U_0,\widetilde \Pi_{\P_0}
\U(t_1)\big) -\Pi_{\P_0} \bt\big(\U_0,\U(t_1)\big) \Big) \\
&+&\Pi_{\P_0} \bt\big(\U_0,\U(t_1)\big) \, ,
\end{eqnarray*}
we have
  $B_3=B_{31}+B_{32}$
with
\begin{equation*}
  B_{31}= \Big( \bt_h\big(\Pi_{\mathbf{RT_0}}\U_0,\widetilde
\Pi_{\P_0} \U(t_1)\big) -\Pi_{\P_0} \bt\big(\U_0,\U(t_1)\big) ,2
\Del \U^1_h-\Del \U^0_h \Big) \, ,
\end{equation*}
\begin{equation*}
B_{32}=\Big( \Pi_{\P_0} \bt\big(\U_0,\U(t_1)\big) ,2 \Del
\U^1_h-\Del \U^0_h \Big).
\end{equation*}
We have
\begin{equation*}
  |B_{31}|   \le \| \bt_h\big(\Pi_{\mathbf{RT_0}}\U_0,\widetilde
\Pi_{\P_0} \U(t_1)\big) -\Pi_{\P_0} \bt\big(\U_0,\U(t_1)\big)
\|_{-1,h} \, \|2 \Del \U^1_h-\Del \U^0_h\|_h
\end{equation*}
So that, using proposition \ref{eq:consbh}
 $|B_{31}|  \le C \, h \, \|\U_0\|_2 \, \|\U(t_1)\|_2 \, \|2 \,
\Del \U^1_h-\Del \U^0_h\|_h$.
Using proposition \ref{propinvp0}   we obtain
\begin{equation*}
 |B_{31}|
\le C  \, \|\U_0\|_2 \, \|\U\|_{{\mathcal{C}}(0,T;\H^2)} \, |2 \,
\Del \U^1_h-\Del \U^0_h|.
\end{equation*}
Let us now bound $B_{32}$. Using the Cauchy-Schwarz inequality and
the stability of $\Pi_{\P_0}$ for the  $\L^2$ norm, we have
\begin{equation*}
|B_{32}| \le \left|\Pi_{\P_0} \bt\big(\U_0,\U(t_1)\big) \right| \,
|2 \, \Del \U^1_h-\Del\U^0_h|  \le \left|\bt\big(\U_0,\U(t_1)\big)
\right| \, |2 \, \Del \U^1_h-\Del\U^0_h|.
\end{equation*}
Integrating by parts, we deduce from   (\ref{eq:defbt})
\begin{equation*}
 \left|\bt\big(\U_0,\U(t_1)\big) \right| \le \sum_{i=1}^2 \left| \U_0 \cdot \nabla u_i(t_1) \right| \le
 |\U_0| \, \|\U(t_1)\|_2 \le C \, |\U_0| \, \|\U\|_{{\mathcal{C}}(0,T;\H^2)}.
\end{equation*}
Thus
 \label{eq:estiv18}
 $|B_{32}| \le C \, |2 \, \Del
\U^1_h -\Del \U^0_h|$.
By gathering the estimates for $B_1$, $B_2$, $B_3$  we get
\begin{equation*}
\label{eq:estiv19}
 \left|\b_h(\U^0_h,\U^1_h,2 \, \Del \U^1_h -\Del \U^0_h)\right| \le C.
\end{equation*}
Thus we have bounded the right-hand side in (\ref{eq:estiv105}). We
 infer from it
\begin{equation*}
   |2 \, \Del \U^1_h-\Del \U^0_h| \le C \, k.
\end{equation*}
Plugging this  estimate into (\ref{eq:estiv104}) and using
(\ref{eq:estincv}) for $m=1$, we get
 (\ref{eq:estincv}) for $m=2$.
Let $m\in\{3,\dots,N-1\}$. We assume that the induction hypothesis
is satisfied up to  rank $n=m-1$.
Let us substract equation (\ref{eq:mombdf}) with the same for $n-1$.
Since the operator $\bt_h$ is bilinear we get
\begin{eqnarray*}
\label{eq:estiv12}
 \frac{3 \, \Del \Ut^{n+1}_h -4 \, \Del \U^n_h+ \Del \U^{n-1}_h}{2 \, k} &-& \frac{1}{\hbox{Re}} \, \Deltat_h
(\Del \Ut^{n+1}_h)+\bt_h(2 \, \Del \U^n_h-\Del
\U^{n-1}_h,\Ut^{n+1}_h)
  \nonumber \\
&+&\bt_h(2 \, \U^n_h-\U^{n-1}_h,\Del \Ut^{n+1}_h)  + \nabla_h
(\delta p^n_h) =\Del \F^{n+1}_h.
\end{eqnarray*}
Let us take the scalar product with  $4 \, k\, \Del \Ut^{n+1}_h$. We
get
\begin{align*}
\label{eq:estiv13}  \left( \frac{3 \, \Del \Ut^{n+1}_h -4 \, \Del
\U^n_h+ \Del \U^{n-1}_h}{2 \, k},4 \, k \,
 \Del \Ut^{n+1}_h \right) &- \frac{4 \, k}{\hbox{Re}} \, \big(\Deltat_h (\Del \Ut^{n+1}_h),\Del \Ut^{n+1}_h\big)
  \nonumber \\
 + 4 \, k \, \b_h(2 \, \Del \U^n_h-\Del \U^{n-1}_h,\Ut^{n+1}_h,\Del
\Ut^{n+1}_h)
&+4 \, k \, \b_h(2 \, \U^n_h-\U^{n-1}_h,\Del \Ut^{n+1}_h,\Del \Ut^{n+1}_h) \\
+4 \, k \, \big(\nabla_h (\delta p^n_h),\Del \Ut^{n+1}_h \big) &= 4
\, k \, (\Del \F^{n+1}_h,\Del \Ut^{n+1}_h).
\end{align*}
According to  proposition \ref{prop:stabbth} we have
\begin{equation*}
\left|4 \, k \, \b_h(2 \, \Del \U^n_h-\Del
\U^{n-1}_h,\Ut^{n+1}_h,\Del \Ut^{n+1}_h)\right| \le C \, k \, |2 \,
\Del \U^n_h-\Del \U^{n-1}_h| \, \|\Ut^{n+1}_h\|_h \, \|\Del
\Ut^{n+1}_h\|_h.
\end{equation*}
Using the  induction hypothesis we get
\begin{equation*}
\left|4 \, k \, \b_h(2 \, \Del \U^n_h-\Del
\U^{n-1}_h,\Ut^{n+1}_h,\Del \Ut^{n+1}_h)\right| \le C \, k ^2 \,
\|\Ut^{n+1}_h\|_h \, \|\Del \Ut^{n+1}_h\|_h.
\end{equation*}
Using the  Young inequality and (\ref{eq:estl2h1}) we infer that
\begin{equation*}
\left|4 \, k \, \b_h(2 \, \Del \U^n_h-\Del
\U^{n-1}_h,\Ut^{n+1}_h,\Del \Ut^{n+1}_h)\right|
 \le \frac{k}{\hbox{Re}} \, \|\Del \Ut^{n+1}_h\|^2_h
+C \, k^2.
\end{equation*}
The other terms are treated like the case  $m=2$. We finally obtain
(\ref{eq:estincv}).
 \qed
\begin{theo}
\label{theo:estp} We assume that the initial values of the scheme
fulfull {\bf (H2)}. There exists a  constant $C>0$ such that for all
$m\in\{2,\dots,N\}$
\begin{equation*}
\label{eq:eqestp}
  k\sum_{n=2}^m |\Pi_{P^{nc}_1} p^n_h|^2 \le C.
  \end{equation*}
\end{theo}
\noindent {\sc Proof.} Let $m \in \{2,\dots,N\}$. We set  $n=m-1$.
Using the  inf-sup condition (\ref{eq:infsd}) and proposition
\ref{prop:propadjh}, we get that there exists $\V_h \in
\P_0\backslash\{\mathbf{0}\}$ such that
\begin{equation}
\label{eq:stabpbdf456}
 C \, \|\V_h\|_h \, |\Pi_{P^{nc}_1} p^{n+1}_h|
 \le -(p^{n+1}_h,\hbox{\rm div}_h \, \V_h)=(\nabla_h p^{n+1}_h,\V_h).
\end{equation}
Plugging (\ref{eq:projib}) into  (\ref{eq:mombdf}) we have
\begin{equation*}
\label{eq:stabpbdf455}
 \nabla_h p^{n+1}_h   = -\frac{3 \, \U^{n+1}_h - 4 \, \U^n_h+\U^{n-1}_h}{2 \, k}
+\frac{1}{\hbox{Re}} \, \Deltat_h \Ut^{n+1}_h -\bt_h(2 \, \U^n_h -
\U^{n-1}_h,\Ut^{n+1}_h)+ \F^{n+1}_h.
\end{equation*}
so that
\begin{align*}
 (\nabla_h p^{n+1}_h,\V_h)
 &= -\left(\frac{3 \, \U^{n+1}_h - 4 \, \U^n_h+\U^{n-1}_h}{2 \, k},\V_h\right)
+\frac{1}{\hbox{Re}}\left(\Deltat_h \Ut^{n+1}_h,\V_h\right) \\
 &-\b_h(2 \, \U^n_h-\U^{n-1}_h,\Ut^{n+1}_h,\V_h)
 + (\F^{n+1}_h, \V_h).
\end{align*}
Using the Cauchy-Schwarz inequality, (\ref{eq:inpoinp0}) and
(\ref{eq:propfh}) we have
\begin{equation*}
\left|\left(\frac{3 \, \U^{n+1}_h - 4 \, \U^n_h+\U^{n-1}_h}{2 \,
k},\V_h\right)\right| \le C \, \left|\frac{3 \, \U^{n+1}_h - 4 \,
\U^n_h+\U^{n-1}_h}{2 \, k}\right| \, \|\V_h\|_h
\end{equation*}
and
\begin{equation*}
  (\F^{n+1}_h,\V_h) \le |\F^{n+1}_h| \, |\V_h| \le C \, |\V_h| \le C  \, \|\V_h\|_h \, ,
\end{equation*}
Thanks to  proposition \ref{prop:stabbth} and theorem
\ref{theo:estv} we have
\begin{equation*}
 \left|\b_h(2 \, \U^n_h-\U^{n-1}_h,\Ut^{n+1}_h,\V_h)\right| \le \left(2 \, |\U^n_h|+|\U^{n-1}_h|\right)
  \, \|\Ut^{n+1}_h\|_h \, \|\V_h\|_h
  \le C \, \|\Ut^{n+1}_h\|_h \, \|\V_h\|_h.
\end{equation*}
And according  to proposition \ref{prop:coerlap} we have
  $\left(\Deltat_h \Ut^{n+1}_h,\V_h\right) \le \|\Ut^{n+1}_h\|_h \, \|\V_h\|_h$.
Thus
\begin{equation*}
\label{eq:stabsbdf}
 (\nabla_h p^{n+1}_h,\V_h)    \le C +C\left( \frac{|3 \, \U^{n+1}_h - 4 \, \U^n_h+\U^{n-1}_h|}{2 \, k}
+\|\Ut^{n+1}_h\|_h\right) \|\V_h\|_h.
\end{equation*}
Comparing with (\ref{eq:stabpbdf456}) we get
\begin{equation*}
 |\Pi_{P^{nc}_1} p^{n+1}_h|       \le C+C \left(
\frac{|3 \, \U^{n+1}_h - 4 \, \U^n_h+\U^{n-1}_h|}{2 \, k}
+\|\Ut^{n+1}_h\|_h  \right).
\end{equation*}
Squaring and  summing from $n=1$ to $m-1$ we obtain
\begin{equation*}
\label{eq:stabpbdf6}
  k\sum_{n=2}^{m} |\Pi_{P^{nc}_1} p^{n}_h|^2
 \le   C+C \, k\sum_{n=1}^{m-1} \frac{|3 \, \U^{n+1}_h - 4 \, \U^n_h+\U^{n-1}_h|^2}{4 \, k^2}+
C \, k\sum_{n=1}^{m-1} \|\Ut^{n+1}_h\|^2_h.
\end{equation*}
The last term on the right-hand side is bounded, thanks to theorem
\ref{theo:estv}. And since
\begin{equation*}
3 \, \U^{n+1}_h - 4 \, \U^n_h+\U^{n-1}_h=3(\U^{n+1}_h - \U^n_h)
-(\U^n_h - \U^{n-1}_h) =3 \, \Del \U^{n+1}_h - \Del \U^n_h
\end{equation*}
we deduce from lemma \ref{lem:estiv}
\begin{equation*}
\label{eq:stabpbdf46}
 \hspace{2.1cm} k\sum_{n=1}^{m-1} \frac{|3 \, \U^{n+1}_h - 4 \, \U^n_h+\U^{n-1}_h|^2}{4 \, k^2}
 \le C \, k\sum_{n=1}^{m} \frac{|\Del \U^n_h|^2}{k^2} \le C. \hspace{2.1cm} \qed
\end{equation*}




\end{document}